\documentclass[letterpaper]{article}
\usepackage[margin=2cm]{geometry}
\usepackage{array}
\usepackage[leqno]{amsmath}
\usepackage{amsfonts,amssymb,amsthm}
\usepackage{amscd,amsxtra}

\usepackage{eucal}
\usepackage{mathrsfs} 

\usepackage[all]{xy}
\SelectTips{cm}{} 

\newtheorem{theorem}{Theorem}[section]
\newtheorem{proposition}[theorem]{Proposition}
\newtheorem{lemma}[theorem]{Lemma}
\newtheorem{corollary}[theorem]{Corollary}
\theoremstyle{definition}
\newtheorem{definition}[theorem]{Definition}

\theoremstyle{remark} \newtheorem{remark}[theorem]{Remark}

\numberwithin{equation}{section}

\newcommand{\curly}[1]{\mathscr{#1}}
\newcommand{\cD}{\curly{D}}
\newcommand{\cL}{\curly{L}}

\newcommand{\field}[1]{\ensuremath{\mathbb{#1}}}
\newcommand{\ZZ}{\field{Z}}
\newcommand{\RR}{\field{R}}
\newcommand{\CC}{\field{C}}
\newcommand{\HH}{\field{H}}
\newcommand{\PP}{\field{P}}
\newcommand{\TT}{\field{T}}

\newcommand{\projspace}[1]{\PP^{#1}}
\newcommand{\pione}{\projspace{1}}

\newcommand{\CR}[4]{[#1\colon #2\colon #3\colon #4]}

\newcommand{\complex}[1]{\mathsf{#1}} 
\newcommand{\CCC}{\complex{C}}

\newcommand{\cat}[1]{\mathsf{#1}}

\newcommand{\sheaf}[1]{\underline{\mathnormal{#1}}}
\newcommand{\sha}[2][\bullet]{\sheaf{A}_{#2}^{#1}}
\newcommand{\shomega}[2][\bullet]{\sheaf{\Omega}_{#2}^{#1}}
\newcommand{\she}[2][\bullet]{\sheaf{\mathcal{E}}_{#2}^{#1}}
\newcommand{\sho}[1]{\mathcal{O}_{#1}}
\newcommand{\deligne}[3][\bullet]{#2(#3)^{#1}_\mathcal{D}}
\newcommand{\deltilde}[3][\bullet]{%
  \smash[t]{\widetilde{#2(#3)}}^{#1}_\mathcal{D}}

\newcommand{\dhh}[2][\bullet]{D(#2)_\mathit{h.h.}^{#1}}
\newcommand{\delH}[4][\bullet]{H^{#1}_\mathcal{D}(#2, #3(#4))}
\newcommand{\dhhH}[3][\bullet]{%
  H^{#1}_{\mathcal{D}_\mathit{h.h.}}(#2,#3)}

\newcommand{\qi}{\xrightarrow{\simeq}}
\newcommand{\lqi}{\overset{\simeq}{\longrightarrow}}

\newcommand{\hyper}[1]{\mathbf{#1}}
\newcommand{\HHH}{\hyper{H}}
\newcommand{\RRR}{\hyper{R}}

\newcommand{\cover}[1]{\mathcal{#1}}

\newcommand{\del}{\partial} 
\newcommand{\delb}{\Bar\partial} 
\newcommand{\deltacheck}{\Check\delta} 

 \DeclareMathOperator{\Tot}{Tot}

\DeclareMathOperator{\cone}{Cone} 
\DeclareMathOperator{\im}{Im} 
\DeclareMathOperator{\PSL}{PSL} 
 
\DeclareMathOperator{\deck}{Deck}

\newcommand{\tame}[2]{\bigl(#1,#2\bigr]}
\newcommand{\tate}{2\pi\sqrt{-1}}

\newcommand{\dual}[2]{\langle#1\,,\,#2\rangle}
\newcommand{\abs}[1]{\lvert#1\rvert}
\newcommand{\norm}[1]{\lVert#1\rVert}

\newcommand{\onehalf}{\frac{1}{2}}
\newcommand{\onefourth}{\frac{1}{4}}
\newcommand{\ihalf}{\frac{\sqrt{-1}}{2}}
\newcommand{\eqdef}{\overset{\mathrm{def}}{=}}
\newcommand{\pic}[1]{\mathrm{Pic}(#1)}
\newcommand{\cm}[1]{\mathscr{CM}(#1)}
\newcommand{\li}{\mathit{Li}_2}
\newcommand{\bwli}{{\cL}_2}
\newcommand{\bwd}{{\cD}_2}

\newcommand{\bei}{Be\u\i{}linson}
\newcommand{\cech}{\v{C}ech}

\title{On hermitian-holomorphic classes related to uniformization, the
  dilogarithm, and the Liouville Action}

\author{Ettore Aldrovandi\\
  Department of Mathematics\\
  Florida State University\\
  Tallahassee, FL 32306-4510, USA\\
  \texttt{aldrovandi@math.fsu.edu} }

\date{}

\begin{document}

\maketitle

\begin{abstract}
  Metrics of constant negative curvature on a compact Riemann surface
  are critical points of the Liouville action functional, which in
  recent constructions is rigorously defined as a class in a
  \cech-de~Rham complex with respect to a suitable covering of the
  surface.
  
  We show that this class is the square of the metrized holomorphic
  tangent bundle in hermitian-holomorphic Deligne cohomology. We achieve
  this by introducing a different version of the hermitian-holomorphic
  Deligne complex which is nevertheless quasi-isomorphic to the one
  introduced by Brylinski in his construction of Quillen line bundles. We
  reprove the relation with the determinant of cohomology construction.
  
  Furthermore, if we specialize the covering to the one provided by a
  Kleinian uniformization (thereby allowing possibly disconnected
  surfaces) the same class can be reinterpreted as the transgression of
  the regulator class expressed by the Bloch-Wigner dilogarithm.
\end{abstract}

\tableofcontents

\section{Introduction}
\label{sec:introduction}

Metrics of constant negative curvature play a very important role in
uniformization problems for compact Riemann surfaces of genus $g>1$. The
condition that the scalar curvature associated to a conformal metric on
a Riemann surface $X$ be equal to $-1$ is equivalent to the fact that
the associated conformal factor satisfies a nonlinear partial
differential equation known as the Liouville equation.

The Liouville equation appears as early as in one of the approaches
considered by Poincar\'e to attack the uniformization theorem
\cite{poincare:liou}. In relatively recent times, it has received
considerable attention in Theoretical and Mathematical Physics due to
the key role it plays in Polyakov's approach to String Theory
\cite{polyakov:bosonic}, especially from the point of view of
non-critical strings and two-dimensional quantum gravity. In this
context one refers to the conformal factor of the metric as the
Liouville ``field.''

As usual in the context of differential equations with a physical
motivation, one would normally like to formulate a variational principle
to express the Liouville equation as an extremum condition. Namely,
given a Riemann surface $X$ and the space $\cm{X}$ of all conformal
metrics on it, the metric of constant negative curvature should be a
critical point of a functional defined over $\cm{X}$. This functional is
the Liouville action. As it happens, action functionals may turn out to
be even more relevant than the equations they are associated to. The
Liouville action is no exception in this sense: it has deep connection
with the geometry of Teichm\"uller spaces
\cite{zogtak1987-1,zogtak1987-2}, and in Physics it describes the
conformal anomaly in String Theory.

Providing a rigorous mathematical definition of the Liouville action
functional is however far from trivial. The very geometric properties of
the Liouville field itself prevent expressing the corresponding
functional as a plain integral of a $2$-form on a Riemann surface.
Correction terms are required, typically in the form of integration of
lower degree forms over the $1$-skeleton of an appropriate simplicial
realization of $X$. (One should notice that this behavior is not
specific to the Liouville equation, and it is by now possible to give a
characterization, in terms of homological algebra, of these type of
functionals, see ref.\ \cite{MR1908413}.)

It is possible to directly determine the necessary correction terms by
requiring that the variational problem be well defined. This, however,
is not completely satisfactory from the point of view of certain
applications to deformation theory, where a consistent definition across
a \emph{family} of surfaces is required. Quite recently, a more
systematic construction, based on the homological algebra techniques
developed by the author and L.\ A.\ Takhtajan in \cite{aldtak1997}, was
carried out by L.\ A.\ Takhtajan and L.-P.\ Teo in ref.\ 
\cite{math.CV/0204318}, generalizing the earlier results of
\cite{zogtak1987-1,zogtak1987-2}.  The authors of ref.\ 
\cite{math.CV/0204318} constructed a \cech\ cocycle with respect to the
\'etale cover of $X$ associated to a quasi-Fuchsian (and more generally
Kleinian) uniformization. Since their construction works across
(Kleinian) deformations, it could be exploited to obtain results of
global nature on the analytic geometry of Kleinian deformation spaces.
As a further result, the authors of loc.\ cit.\ were able to rigorously
prove the validity of the ``holography principle'' for the Liouville
action corresponding to a large class of Kleinian (in particular
Fuchsian and quasi-Fuchsian) uniformizations.  Specifically, they proved
that given a second kind Kleinian group, the corresponding Liouville
action can be obtained as the regularized limit of the hyperbolic volume
of the corresponding associated $3$-manifold.  This extends to the
general Kleinian case a previous formula obtained by Krasnov
\cite{MR2002k:81230} for classical Schottky groups.

Our interest in this matter is two-fold. From the perspective of the
newer methods adopted in \cite{aldtak2000}, the covering map $U\to X$
associated to the uniformization by a discrete group $\Gamma=\deck
(U/X)$ is but one of the many possible covers comprising an appropriate
category $\cat{C}$ of, say, local diffeomorphisms $U\to X$---the most
obvious choice being that of standard open cover $\cover{U}=\{ U_i\}$
with associated space $U=\coprod_i U_i$. In particular one expects to be
able to apply the methods of \cite{aldtak1997} and \cite{aldtak2000}
uniformly on a class of reasonably behaved covers of $X$.

Second, the focus of ref.\ \cite{aldtak2000} was on the rigorous
definition of a functional for quasi-conformal deformation of the
Riemann surface $X$ and its application to the study of projective
structures. A main result is that the construction of the action is
possible thanks to the vanishing of the ``tame symbol'' (see refs.\ 
\cite{del:symbole} and \cite{brymcl:deg4_II} for the relevant
definitions) $\tame{T_X}{T_X}$, where $T_X$ is the holomorphic tangent
line bundle of $X$. The vanishing determines local choices (with respect
to a cover) of a Bloch-type dilogarithm which then allow for a
cohomological construction of the action. There are many indication that
the Liouville action ought to be the hermitian square of a functional of
the type studied in \cite{aldtak2000}.\footnote{From a physical point of
  view this originates in the modular geometry approach to Conformal
  Field Theory advocated by Friedan and Shenker in \cite{MR88b:81146}.
  Mathematically speaking, it is one of the many proposed forms of the
  holomorphic factorization property for determinant line bundles.} Thus
it is natural to ask whether there is an analogous mechanism as the one
in loc.\ cit.\ to obtain a general construction of the Liouville action
by replacing the holomorphic symbol maps and dilogarithms with
corresponding real objects.

In this paper we answer this question in the affirmative. More
precisely, we show that the Liouville action (up to the area term
which is given by an ordinary $2$-form) can be computed as a
symbol map taking values in hermitian holomorphic Deligne
cohomology, first introduced by Brylinski and McLaughlin in their
study of degree four characteristic classes
\cite{brymcl:deg4_II}. (By way of comparison, the tame symbols
used in ref.\ \cite{aldtak2000} used holomorphic and smooth
Deligne cohomology.) In particular we show that the dilogarithm
type terms are replaced here by the Bloch-Wigner function, the
real valued counterpart of the dilogarithm (see refs.
\cite{MR2001i:11082}, and \cite{MR94k:19002,MR2002g:52013} for a
review.)

The appearance of the Bloch-Wigner function ties very well with
the holography property of the Liouville function proved in
\cite{math.CV/0204318} in the following sense. As mentioned
before, the Liouville action (up to the area term) relative to a
Kleinian uniformization\footnote{Note that $X$ is allowed to be
  disconnected.}  $U\to X$ can also be computed as the
``regularized volume'' of the associated $3$-manifold $N=\Gamma
\backslash (U \cup \HH^3)$, where $\Gamma =\deck (U/X)$ as
before, $U\subset \PP^1$ is the domain of discontinuity for
$\Gamma$, and $\HH^3$ is the standard hyperbolic $3$-space. (To
define the regularized volume would lead us too far afield. It
suffices to mention that the conformal factor of a metric on
$X=\partial N$ can be used to select a compact submanifold
$N_\epsilon$ whose volume is finite. One then subtracts from the
volume of $N_\epsilon$ the areas of the boundary components and
other carefully chosen constants independent of the metric
structure, so that the resulting quantity will have a finite
limit as $\epsilon \to 0$.)  On the other hand, the hyperbolic
volume in three dimensions corresponds to a three dimensional
(purely imaginary) class on $\PSL_2(\CC)$ expressible through the
Bloch-Wigner dilogarithm, the so-called \emph{regulator class.}
We show that the regulator is precisely the class that needs to
be killed in order to close the cohomological descent conditions
required to calculate the Liouville action for a covering map
$U\to X$ with covering group a Kleinian group $\Gamma$.  This is
possible, since for a second kind Kleinian group the quotient
$\HH^3/\Gamma$ is non-compact, hence it carries no cohomology in
dimension three, so the class represented by the Bloch-Wigner
function, pulled back to $\Gamma$ via the imbedding $\Gamma
\hookrightarrow \PSL_2(\CC)$, vanishes.

Returning to the cohomological interpretation of the construction
of the Liouville action, it should also be noted that leaving
aside the area term, our results show that the cohomologically
non trivial part is indeed a square. Namely, for a conformal
metric $\rho \in \cm{X}$ we consider the pair $(T_X,\rho)$ as a
holomorphic line bundle equipped with an hermitian metric. Then,
using that hermitian holomorphic Deligne cohomology has a cup
product, we show that the Liouville action is just the square of
the class of $(T_X,\rho)$. In fact this identification holds at
the level of cocycles, rather than only for the corresponding
classes.

Again leaving aside the area term, it immediately follows from
the properties of hermitian holomorphic Deligne cohomology that
most of the story carries over to the case of a \emph{pair} of
holomorphic line bundles $L$ and $L'$ equipped with hermitian
metrics $\rho$ and $\rho'$, respectively. Furthermore, Brylinski
shows in \cite{bry:quillen} that the pairing of two such
holomorphic line bundles with metrics corresponds to the pairing
defined by Deligne on the determinant of cohomology in
\cite{MR89b:32038}. Without introducing the machinery of
$2$-gerbes, we reobtain this result in our setting. Specifically,
we directly obtain Gabber's formula for the hermitian metric on
the determinant line from the explicit cocycle for the cup
product of two metrized line bundles. In turn this shows that the
Liouville action is a multiple of the determinant of cohomology,
thereby generalizing earlier results (cf.\ 
\cite{zograf1990})---without assuming criticality.

\subsection{Organization of the paper}
\label{sec:Organization-paper}

This paper is organized as follows.
Sections~\ref{sec:deligne-complexes}
and~\ref{sec:herm-holom-deligne} are devoted to expounding some
background material for the sake of keeping this paper
self-contained and to put the reader in position of reproducing
the necessary calculations. Section~\ref{sec:deligne-complexes}
contains background facts on Deligne cohomology, paying special
attention to the product structures and the cone constructions.
We provide some examples and collect some facts about the
dilogarithm from the point of view of Deligne cohomology.  The
particular model of hermitian holomorphic Deligne cohomology we
use later in the paper requires certain constructions available
in the literature, and recalled in
section~\ref{sec:deligne-complexes}, to be slightly modified in
order to obtain a (graded) commutative product. The necessary
arguments, being somewhat outside the line of development of the
paper are presented in Appendix~\ref{sec:Cones}. Hermitian
holomorphic Deligne cohomology is introduced in
section~\ref{sec:herm-holom-deligne}. We give the definition as
in refs.\ \cite{brymcl:deg4_II} and \cite{bry:quillen}, and then
introduce another model which, albeit more complex, has the
advantage for us of keeping the metric structure explicit. We
have explicitly proved the isomorphism in
Lemma~\ref{lemma:qi-dhh-bry}. For the sake of completeness, we
give an explicit description of the cocycle determined by a
holomorphic line bundle with hermitian metric, and in
sect.~\ref{sec:Cup-prod-herm} we explicitly compute the cup
product of two metrized line bundles for later usage. Results
about the existence of a fiber integration map are mentioned in
the paper, so some background material is provided in
sect.\ \ref{sec:integr-along-fiber}.

Sections~\ref{sec:Conf-metr-Liouv} and~\ref{sec:conf-metr-herm}
form the core of the paper. The direct construction of the
Liouville action according to the techniques of refs.\ 
\cite{aldtak1997,aldtak2000,math.CV/0204318} is presented in
section~\ref{sec:Conf-metr-Liouv}. Since explicit calculations
have been presented in great detail in ref.\ 
\cite{math.CV/0204318}, and the calculations we need are quite
straightforward, we keep details to a minimum. In
Proposition~\ref{prop:2} and Corollary~\ref{cor:1} we show that
the Liouville action functional computed via descent theory does
solve the variational problem. These results have appeared also
in ref.\ \cite{math.CV/0204318} and are presented for here
completeness, although our choice for the cover is different. The
framework of loc.\ cit.\ is that of a Kleinian cover $U\to X\cong
U/\Gamma$, where $\Gamma$ is a second kind purely loxodromic
geometrically finite Kleinian group, which then we treat in some
detail in sect.\ \ref{sec:liouv-funct-bloch}. We illustrate how
the genuine Bloch-Wigner function of sect.\ 
\ref{sec:remarks-dilogarithm} appears in the descent equations
relative to the Kleinian cover. Having observed descent equations
close on general cohomological grounds in sect.\ 
\ref{sec:solut-vari-probl}, we now point out that for the case of
a Kleinian cover this is due to the vanishing of the regulator
class for the non-compact $3$-manifold $\HH^3/\Gamma$.

A true geometric construction of the Liouville action, which does
not rely on the arguments of sect.\ \ref{sec:Conf-metr-Liouv} to
close the descent equations, is carried out in
section~\ref{sec:conf-metr-herm}, which contains our main result:
We compare the descent calculations with the cup products
computed in section~\ref{sec:herm-holom-deligne} and conclude
that the quadratic part of the Liouville action is in fact
(modulo an area term) the cup square of the metrized holomorphic
tangent bundle in Hermitian holomorphic Deligne cohomology, see
Theorem~\ref{thm:2}, Corollary~\ref{corollary:1}, and
Proposition~\ref{prop:3}. It follows that from this point of view
the descent equations close automatically, without the need for
special arguments.  We then prove that the cup square is
identified with the determinant of cohomology construction in
Theorem~\ref{thm:3} and Corollary~\ref{corollary:2}. Auxiliary
facts about the homological algebra of cones and conventions on
Kleinian groups are stored in the appendices.

Finally, we draw our conclusions in sect.\ 
\ref{sec:conclusions-outlook}.

\subsection{Notation and conventions}
\label{sec:notation}

If $z$ is a complex number, then $\pi\sb{p} (z) \eqdef \onehalf ( z +
(-1)\sp p \Bar z)$, and similarly for any other complex quantity, e.g.
complex valued differential forms. If $A$ is a subring of $\RR$, we will
use the ``twist'' $A(j) = (2\pi \sqrt{-1})^j\,A$.

If $X$ is a complex manifold, $\sha{X}$ and $\shomega{X}$ denote
the de~Rham complexes of smooth $\CC$-valued and holomorphic
forms, respectively. We denote by $\she{X}$ the de~Rham complex
of sheaves of real valued differential forms and by $\she{X} (j)$
the twist $\she{X} \otimes_\RR \RR(j)$. We set $\sho{X} \equiv
\shomega[0]{X}$ as usual.  When needed, $\sha[{p,q}]{X}$ will
denote the sheaf of smooth $(p,q)$-forms. We use the standard
decomposition $d=\del + \delb$ according to types. Furthermore,
we introduce the differential operator $d^c = \del -\delb$
(contrary to the convention, see, e.g. \cite{lang:arakelov}, we
omit the factor $1/(4\pi \sqrt{-1})$). We have $2\del\delb =
d^cd$. The operator $d^c$ is an imaginary one, and accordingly,
we have the rules
\begin{gather*}
  d\pi_p(\omega) = \pi_p(d\omega)\\
  d^c\pi_p(\omega) = \pi_{p+1}(d^c\omega)
\end{gather*}
for any complex form $\omega$.

An open cover of $X$ will be denoted by $\cover{U}_X$. If $\{U_i\}_{i\in
  I}$ is the corresponding collection of open sets, we denote $U_{ij} =
U_i\cap U_j$, $U_{ijk} = U_i\cap U_j\cap U_k$, and so on. We can also
consider more general covers $\cover{U}_X = \{ U_i \to X\}_{i\in I}$
where the maps are regular coverings in the appropriate category. In
this case intersections are replaced by $(n+1)$-fold fibered products
\begin{math}
  U_{i_0 i_1\dots i_n}= U_{i_0} \times_X\dots \times_X U_{i_n}\,.
\end{math}
Open coverings fit this more general description, since if $U_i$ and
$U_j$ are two open sets, then $U_i\cap U_j = U_i\times_X U_j$, where the
fiber product is taken with respect to the inclusion maps. As another
example, one can consider regular covering maps $U\to X$ with $\Gamma =
\mathrm{Deck}(U/X)$, and in this case, taking the fiber product over $X$
$(n+1)$-times, one gets
\begin{math}
  U\times_X\dots\times_X U = U\times \Gamma \times \dots \times \Gamma\,,
\end{math}
where the group factor is repeated $n$-times. This includes the cases of
Kleinian (and in particular Fuchsian) covers of Riemann surfaces.

The \emph{nerve} of the cover $\cover{U}_X$ is
the simplicial object
\begin{math}
  n\mapsto N_n(\cover{U}_X) = \coprod U_{i_0} \times_X\dots \times_X
  U_{i_n}
\end{math}
where $N_n(\cover{U}_X)$ maps into $N_{n-1}(\cover{U}_X)$ in $(n+1)$
ways by forgetting in turn each factor. For open covers this just yields
the expected inclusion maps.

If $\sheaf{F}^\bullet$ is a complex of sheaves on $X$, its \cech\ 
resolution with respect to a covering $\cover{U}_X\to X$ is the double
complex
\begin{math}
  \CCC^{p,q} (\sheaf{F}) \eqdef \check{C}^q(\cover{U}_X,\sheaf{F}^p)\,,
\end{math}
where the $q$-cochains with values in $\sheaf{F}^p$ are given by
\begin{math}
\prod \sheaf{F}^p (U_{i_0} \times_X\dots \times_X U_{i_n})\,.
\end{math}
The \cech\ coboundary operator is denoted $\deltacheck$. The sign
convention we are going to use is that the index along the \cech\ 
resolution is the \emph{second} one, so if we denote by $d$ the
differential in the complex $\sheaf{F}^\bullet$, the total differential
in the total simple complex of $\check{C}^q(\cover{U}_X,\sheaf{F}^p)$
will be $D=d\pm \deltacheck$. For open covers we just get the familiar
\cech\ (hyper)cohomology. The other interesting example is that of a
regular covering map $U\to X$: \cech\ cohomology with respect to this
cover is the same as group cohomology for $\Gamma = \mathrm{Deck}(U/X)$
with coefficients in the $\Gamma$-module $\sheaf{F}^p(U)$.

The Koszul sign rule that results in a sign being picked whenever two
degree indices are formally exchanged is applied. In particular, for
\cech\ resolutions of complexes of sheaves, it leads to the following
conventions. If $\sheaf{G}^\bullet$ is a second complex of sheaves on
$X$, then one defines the cup product
\begin{equation*}
  \cup : \CCC^{p,q}(\sheaf{F}) \otimes \CCC^{r,s}(\sheaf{G})
  \longrightarrow \Check{C}^{q+s}(\cover{U}_X,\sheaf{F}^p\otimes
  \sheaf{G}^r) \subset \CCC^{p+r,q+s}(\sheaf{F}\otimes\sheaf{G})
\end{equation*}
of two elements $\{f_{i_0,\dots,i_q}\}\in \CCC^{p,q}(\sheaf{F})$
and $\{g_{j_0,\dots,j_s}\} \in \CCC^{r,s}(\sheaf{G})$ by
\begin{equation*}
  (-1)^{qr}\,f_{i_0,\dots,i_q}\otimes
  g_{i_q,i_{q+1},\dots,i_{q+s}} \,.
\end{equation*}

\section*{Acknowledgments}
Parts of this work were completed during visits at the
International School for Advanced Studies (SISSA) in Trieste,
Italy, and at the Department of Mathematics, Instituto Superior
T\'ecnico in Lisbon, Portugal. I would like to thank both
institutions for support and for creating a friendly and
stimulating research environment. I would also like to thank
Paolo Aluffi, Phil Bowers, Ugo Bruzzo, Johan Dupont, Leon
Takhtajan for illuminating discussions and/or patiently answering
my many questions. Also, special thanks are due to the referee
for his or her thoroughness and for providing very detailed
comments.

\section{Deligne complexes}
\label{sec:deligne-complexes}

\subsection{Cup products on cones}
\label{sec:cup-products-cones}

Recall that the cone of a map $f: X^\bullet \to Y^\bullet$
between two complexes is the complex $C^\bullet (f) =
X^\bullet[1]\oplus Y^\bullet$ with differential $d_f(x,y) =
(-d\,x, f(x) + d\,y)$, where $[k]$ denotes the shift functor. The
cone fits into the exact sequence
\begin{equation*}
  0 \longrightarrow
  Y^\bullet \longrightarrow C^\bullet (f) \longrightarrow
  X^\bullet[1] \longrightarrow 0\,.
\end{equation*}
The following constructions are a special case of those considered by
\bei\ in ref.\ \cite{bei:hodge_coho}. Suppose we are given complexes
$X\sp\bullet\sb i$, $Y\sp\bullet\sb i$, and $Z\sp\bullet\sb i$ and maps
$f\sb i : X\sp\bullet\sb i \rightarrow Z\sp\bullet \sb i$, $g\sb i :
Y\sp\bullet\sb i \rightarrow Z\sp\bullet \sb i$, for $i=1,2,3$.  Suppose
also that we have product maps $X^\bullet_1\otimes X^\bullet_2
\xrightarrow{\cup} X^\bullet_3$, and similarly for $Y^\bullet_i$, and
$Z^\bullet_i$, strictly compatible with the $f_i$, $g_i$ in the obvious
sense. Then we can consider the cones
\begin{equation*}
  \cone \big( f_i -g_i \big)[-1] \equiv
  \cone \big( X^\bullet_i\oplus Y^\bullet_i \xrightarrow{f_i -
  g_i} Z^\bullet_i \big)[-1]\,.
\end{equation*}
For a real parameter $\alpha$, there is a family of products
\begin{equation}
  \label{eq:1}
  \cone \big( f_1 - g_1 \big)[-1] \otimes
  \cone \big( f_2 - g_2 \big)[-1] \xrightarrow{\cup_\alpha}
  \cone \big( f_3 - g_3 \big)[-1]
\end{equation}
determined as follows. For
\begin{math}
  (x_i,y_i,z_i) \in X^\bullet_i \oplus Y^\bullet_i \oplus
  Z^{\bullet -1}_i\,,\;i=1,2\,,
\end{math}
one defines
\begin{equation}
  \label{eq:2}
  \begin{split}
    (x_1,y_1,z_1) \cup_\alpha (x_2,y_2,z_2) =
    \Big(&x_1\cup x_2, y_1\cup y_2, \\
    &(-1)^{\deg (x_1)}
    \big((1-\alpha )f_1(x_1) + \alpha g_1(y_1) \big) \cup z_2 \\
    &\quad +z_1\cup \big( \alpha f_2(x_2) +
    (1-\alpha)g_2(y_2)\big) \Big)\,.
  \end{split}
\end{equation}
Note that $\deg (x_1)=\deg (x_1,y_1,z_1)$. Checking that
$\cup_\alpha$ is a map of complexes is a straightforward routine
calculation. Different products for different values $\alpha$ and
$\beta$ of the real parameter are homotopic. Explicitly, we have
\begin{equation*}
  (x_1,y_1,z_1) \cup_\alpha (x_2,y_2,z_2) -
  (x_1,y_1,z_1) \cup_\beta (x_2,y_2,z_2)  =
  \big( d \, h_{\alpha,\beta} + h_{\alpha,\beta} d \big)
  ((x_1,y_1,z_1) \otimes (x_2,y_2,z_2)) \,,
\end{equation*}
where the homotopy
\begin{equation*}
  h_{\alpha,\beta} : \Tot \Big( \cone \big( f_1 - g_1\big)[-1]
  \otimes \cone \big( f_2 - g_2\big)[-1] \Big)^\bullet
  \longrightarrow \cone \big( f_3 - g_3\big)[-1]^{\bullet -1}
\end{equation*}
is given by the formula
\begin{equation}
  \label{eq:3}
  h_{\alpha,\beta} ((x_1,y_1,z_1) \otimes (x_2,y_2,z_2)) =
  (\alpha -\beta) \, (-1)^{\deg (x_1) -1} (0,0,z_1\cup
  z_2) \,.
\end{equation}
If the products
\begin{math}
  X^\bullet_1 \otimes X^\bullet_2
  \xrightarrow{\cup} X^\bullet_3\,,
\end{math}
etc., are graded commutative, then the swap functor on the tensor
product maps the $\cup_\alpha$ product structure on the cones into the
$\cup_{1-\alpha}$ structure. Using the homotopy~\eqref{eq:3} it follows
at once that there is a well defined graded commutative product in
cohomology.

If we do not assume the product structures $X^\bullet_1\otimes
X^\bullet_2 \xrightarrow{\cup} X^\bullet_3$, etc., are strictly
compatible with the maps $f_i$, $g_i$, some of the preceding
constructions must be slightly modified. With an eye toward
certain constructions to be carried out later in this paper, let
us assume we have compatibility up to homotopy, namely there
exist maps
\begin{gather*}
  h \colon \bigl( X_1\otimes X_2 \bigr)^\bullet
  \longrightarrow Z_3^{\bullet -1}\\
  k \colon \bigl( Y_1\otimes Y_2 \bigr)^\bullet
  \longrightarrow Z_3^{\bullet -1}
\end{gather*}
such that 
\begin{equation}
  \label{eq:4}
  \begin{gathered}
    f_3\circ \cup - \cup \circ (f_1\otimes f_2)
    = d\, h + h\,d \\
    g_3\circ \cup - \cup \circ (g_1\otimes g_2)
    = d\, k + k\,d\,,
  \end{gathered}
\end{equation}
with obvious meaning of the symbols.
\begin{lemma}
  \label{lemma:mod-cup-prod}
  Let $X_i$, $Y_i$, $Z_i$ and the maps $f_i$, $g_i$ be as above. Let
  $\alpha$ be a real parameter. We have a product of type~\eqref{eq:1}
  for the cones $\cone(f_i-g_i)[-1]$ defined by the following
  modification of formula~\eqref{eq:2}:
  \begin{equation}
    \label{eq:5}
    \begin{split}
      (x_1,y_1,z_1) \cup_\alpha (x_2,y_2,z_2) =
      \Big(&x_1\cup x_2, y_1\cup y_2, \\
      &(-1)^{\deg (x_1)}
      \big((1-\alpha )f_1(x_1) + \alpha g_1(y_1) \big) \cup z_2 \\
      &\quad +z_1\cup \big( \alpha f_2(x_2) +
      (1-\alpha)g_2(y_2)\big)\\
      &\qquad -h(x_1\otimes x_2) +k(y_1\otimes y_2)
      \Big)\,.
    \end{split}   
  \end{equation}
  The product~\eqref{eq:5} is a map of complexes and two products
  $\cup_\alpha$ and $\cup_\beta$ are related by the same homotopy
  formula~\eqref{eq:3}.
\end{lemma}
\begin{proof}
  Direct verification.
\end{proof}
This modified framework carries over to the full structure considered by 
\bei\ in ref.\ \cite{bei:hodge_coho}. We will still refer to this
modified product as the \bei\ product. It is also necessary to relax the 
assumption that the products
\begin{math}
  X^\bullet_1 \otimes X^\bullet_2
  \xrightarrow{\cup} X^\bullet_3\,,
\end{math}
etc., be graded commutative. It is possible to complete all the diagrams
so that the permutation of factors in the tensor products still yields a
homotopy commutative product~\eqref{eq:5} for the cones. As a
consequence the induced product in cohomology will still be graded
commutative. Explicit formulas are not needed except to ensure this
latter fact, therefore we shall not discuss this matter any further and
refer the reader to the appendix, where a brief but explicit treatment
can be found.

\subsection{Deligne complexes}
\label{sec:deligne-compl}

Let $X$ be a complex manifold. Recall the standard Hodge
filtration of $\shomega{X}$:
\begin{equation}
  \label{eq:6}
  F^p\shomega{X} : 0 \longrightarrow \cdots \longrightarrow
  \shomega[p]{X} \longrightarrow \cdots \longrightarrow
  \shomega[n]{X}\,,
\end{equation}
where $n=\dim_\CC X$.

The corresponding filtration for the complex of smooth
$\CC$-valued forms is defined as follows: denote by $F^p\sha{X}$
the subcomplex of $\sha{X}$ comprising forms of type $(r,s)$
where $r$ is at least $p$, so that $F^p\sha[k]{X} = \oplus_{r\geq
  p} \sha[{r,k-r}]{X}$. Then (cf. \cite{del:hodge_II}) the
inclusion $\shomega{X} \hookrightarrow \sha{X}$ is a
quasi-isomorphism respecting the filtrations, namely
$F^p\shomega{X} \hookrightarrow F^p\sha{X}$, and the latter
inclusion induces an isomorphism in cohomology.

If $A$ is a subring of $\RR$, and $\imath$ and $\jmath$ denote
the inclusions of $A(p)$ and $F^p\shomega{X}$ into $\shomega{X}$
respectively, the $p$-th Deligne complex of sheaves is
defined by
\begin{equation}
  \label{eq:7}
  \deligne{A}{p} = \cone \bigl( A(p)_X \oplus F^p\shomega{X}
  \xrightarrow{\imath -\jmath} \shomega{X} \bigr)[-1]\,.
\end{equation}
It is quasi-isomorphic to the complex:
\begin{equation}
  \label{eq:8}
  \cone \bigl( A(p)_X \oplus F^p\sha{X}
  \xrightarrow{\imath -\jmath} \sha{X} \bigr)[-1]\,,
\end{equation}
where $\imath$ and $\jmath$ have the same meaning. We also notice
the quasi-isomorphism
\begin{equation}
  \label{eq:9}
  \deligne{A}{p}\lqi
  \bigl( A(p)_X 
  \overset{\imath}{\longrightarrow} \sho{X}
  \overset{d}{\longrightarrow} \shomega[1]{X}
  \overset{d}{\longrightarrow} \dotsm
  \overset{d}{\longrightarrow} \shomega[{p-1}]{X}\bigr)\,.
\end{equation}
When $A=\RR$ there are further quasi-isomorphisms, namely
\begin{equation*}
  \deligne{\RR}{p} \xrightarrow{\simeq} \cone \big(
  F^p\shomega{X} \rightarrow \she{X} (p-1)\big) [-1]
  \xrightarrow{\simeq} \cone \big( F^p\sha{X} \rightarrow \she{X}
  (p-1)\big) [-1]
\end{equation*}
since the maps
\begin{equation*}
  \bigl( \RR(p) \rightarrow \shomega{X}\bigr)
  \xrightarrow{\simeq} \bigl( \RR(p)
  \rightarrow \CC \bigr) \xrightarrow{\simeq} \RR(p-1)
  \xrightarrow{\simeq} \she{X}(p-1)
\end{equation*}
are all quasi-isomorphisms in the derived category, cf.
\cite{esn-vie:del}. Here we have used $\CC \cong \RR(p) \oplus
\RR(p-1)$. Following op.~cit., we set:
\begin{equation}
  \label{eq:11}
  \deltilde{\RR}{p} \eqdef
  \cone \big( F^p\sha{X} \xrightarrow{-\pi\sb{p-1}} \she{X} (p-1)\big)
  [-1] \, .  
\end{equation}
Again, there is an explicit quasi-isomorphism
(\cite{MR86h:11103,esn-vie:del}):
\begin{equation}
  \label{eq:12}
  \begin{gathered}
    \rho_p : \deligne{\RR}{p} \overset{\simeq}{\longrightarrow}
    \deltilde{\RR}{p} \\
    \rho_p \lvert_{\RR(p)} = 0\,,\quad \rho_p
    \lvert_{F^p\shomega{X}} = \mathit{incl} \,,\quad \rho_p
    \lvert_{\shomega{X}} = \pi_{p-1}
  \end{gathered}
\end{equation}

The \emph{Deligne cohomology groups} of $X$ with coefficients in
$A(p)$ are the hypercohomology groups
\begin{equation*}
  \delH{X}{A}{p} = \HHH^\bullet (X, \deligne{A}{p})\,. 
\end{equation*}
and clearly, any complex quasi-isomorphic to $\deligne{A}{p}$
would do. In order to perform calculations with these cohomology
groups we shall normally resort to a \cech\ resolution with
respect to an open cover $\cover{U}_X$ of $X$ or an \'etale map
$\cover{U}_X \to X$, e.g. a regular cover with deck group
$\Gamma$.

One of the important properties of Deligne cohomology is the existence
of a graded commutative cup product
\begin{equation}
  \label{eq:13}
  \delH[i]{X}{A}{p} \otimes \delH[j]{X}{A}{q}
  \xrightarrow{\cup} \delH[i+j]{X}{A}{p+q}\,,
\end{equation}
which follows from the existence of the \bei{} cup product at the
level of Deligne complexes whose construction was recalled above.
There are products $A(p) \otimes A(q) \to A(p+q)$ and
$F^p\shomega{X} \otimes F^q\shomega{X} \to F^{p+q}\shomega{X}$,
plus the obvious (wedge) product on $\shomega{X}$, thus it
follows from the cone version~\eqref{eq:7} that the Deligne
complexes come equipped with the \bei{} product, and therefore
the Deligne cohomology groups inherit the graded commutative cup
product~\eqref{eq:13}. The explicit form, that is, the
translation of~\eqref{eq:2} to the case at hand can be found
in~\cite{esn-vie:del}. The explicit form of the cup product for
the complex~\eqref{eq:11} as computed in \cite{MR86h:11103} (see
also \cite{esn-vie:del}) will be needed in the sequel. Let
$(\omega\sb 1, \eta\sb 1)$ be an element of degree $n$ in
$\deltilde{\RR}{p}$---this means that $\omega\sb 1\in
F^p\shomega[n]{X}$ and $\eta\sb 1\in \she[n-1]{X}(p-1)$---and
$(\omega\sb 2, \eta\sb 2)$ any element in $\deltilde{\RR}{q}$.
The product is defined by the formula:
\begin{equation}
  \label{eq:15}
  (\omega_1,\eta_1) \,\Tilde\cup\, (\omega_2,\eta_2) = \bigl(
  \omega\sb 1 \wedge \omega\sb 2 ,
  (-1)^n\, \pi\sb p\omega\sb 1 \wedge \eta\sb 2
  +\eta\sb 1\wedge \pi\sb q \omega\sb 2 \bigr)\,.
\end{equation}
The product $\Tilde\cup$ is a morphism of complexes and (modulo
the quasi-isomorphisms $\rho_p$) is homotopic to the \bei\ 
products $\cup_\alpha$ on the complexes $\deligne{\RR}{p}$.
Specifically, if we denote an element of $\deligne{\RR}{p}$ by
the triple $(r,f,\omega)$, where $r\in \RR(p)_X$, $f\in
F^{p}\!\shomega{X}$, and $\omega \in \shomega{X}$, the homotopy
between $\Tilde\cup \circ (\rho_p\otimes \rho_q)$ and
$\rho_{p+q}\circ \cup_\alpha$ is given by
\begin{equation}
  \label{eq:16}
  \Tilde h ((r,f,\omega)\otimes (r',f',\omega'))
  = (-1)^{\deg \omega} \bigl(0,
  (1-\alpha)\, \pi_p\omega\wedge\pi_{q-1}\omega'
  -\alpha\, \pi_{p-1}\omega\wedge\pi_{q}\omega'
  \bigr)\,.
\end{equation}

\subsection{Examples}
\label{sec:examples}

\subsubsection{}
Let $A=\ZZ$. From~\eqref{eq:9} we have $\deligne{\ZZ}{1} \qi
\sho{X}^\times[-1]$ via the standard exponential sequence, so
that
\begin{math}
  \delH[k]{X}{\ZZ}{1} \cong H^{k -1}(X, \sho{X}^\times)\,.
\end{math}
In particular $\delH[1]{X}{\ZZ}{1}\cong H^0(X,\sho{X}^\times)$,
the global invertibles on $X$, and $\delH[2]{X}{\ZZ}{1} \cong
\pic{X}$, the Picard group of line bundles over $X$.
  
If an open cover $\{U_i\}_{i\in I}$ of $X$ is chosen, the class
of a line bundle $L$ in $\delH[2]{X}{\ZZ}{1}$ can be represented
via a \cech\ resolution by the cocycle $(f_{ij},c_{ijk})$, where
$f_{ij}\in \sho{X}(U_{ij})$, and $c_{ijk}\in \ZZ(1)_X(U_{ijk})$.
Thus the functions $f_{ij}$ should be interpreted as
\emph{logarithms} of the corresponding transition functions for
$L$. Then, the collection $c_{ijk}=(\deltacheck f)_{ijk}$
provides a representative for the first Chern Class $c_1(L)$.
Similarly, an invertible function $f$ would be described by the
collection $f_i$ of its logarithms on each open $U_i$, plus the
``integers'' $m_{ij} = f_j -f_i\in \ZZ(1)$ on each $U_{ij}$.
  
\subsubsection{}
Still using the exponential sequence, $\deligne{\ZZ}{2} \qi
\bigl( \sho{X}^\times \xrightarrow{d\log}
\shomega[1]{X}\bigr)[-1]$. Thus $\delH[2]{X}{\ZZ}{2}$ is the
group of isomorphism classes of holomorphic line bundles with
(holomorphic) connection. Using the (in fact, any) product
$\deligne{\ZZ}{1} \otimes \deligne{\ZZ}{1} \to \deligne{\ZZ}{2}$,
the cup product of two global invertible holomorphic functions
$f$ and $g$ on $X$ yields a line bundle with connection---the
tame symbol---denoted by $\tame{f}{g}$ whose class is in
$\delH[2]{X}{\ZZ}{2}$,
see~\cite{del:symbole,bloch:dilog_lie,rama:reg_hei}. Higher cup
products in this spirit have been studied
in~\cite{brymcl:deg4_II}.
  
\subsubsection{}
If $A=\RR$, we have
\begin{math}
  \delH[2p]{X}{\RR}{p} = H^{2p}(X,\RR(p)) \cap H^{p,p}(X)\,.
\end{math}
The \cech\ representative $(c_{ijk},f_{ij})$ of a class in
$\delH[2]{X}{\ZZ}{1}$ mentioned above maps to the cocycle $(-d{}
f_{ij}, -\abs{f_{ij}})$ under~\eqref{eq:12}. Taking into account
that the $f_{ij}$ are the logarithms of the transition functions,
the corresponding $(1,1)$ class would be given by the associated
canonical connection, see sec.~\ref{sec:line-bun}.
  
\subsubsection{}\label{sec:ex_4}
$\delH[1]{X}{\RR}{1}$ is the group of real valued functions $f$ on $X$
such that there exists a holomorphic one-form $\omega$ such that $\pi\sb
0\omega = df$. In other words it is the group of those smooth functions
$f$ such that $\del f$ is holomorphic, which amounts to say that such an
$f$ itself is harmonic.

\subsection{Remarks on the cup product $f\cup g$}
\label{sec:remarks-dilogarithm}

It is convenient to consider the case of the cup product of two
invertible functions $f$ and $g$ in various complexes in more detail,
and to introduce some related notions we shall need later.

As observed, $\deligne{\ZZ}{1} \qi \sho{X}^\times[-1]$ and an
invertible function $f$ can be considered as an element of
$\delH[1]{X}{\ZZ}{1}$. Therefore, via~\eqref{eq:12}, it induces
$\rho_1(f) \in \delH[1]{X}{\RR}{1}$ represented by $(d\log f,
\log\,\abs{f})$. (Note that $\log \abs{f}$ fits the description
of $\delH[1]{X}{\RR}{1}$ in \ref{sec:ex_4}.) Given two such $f$
and $g$, the expression for the cup product~\eqref{eq:15} gives
the following element of $\delH[2]{X}{\RR}{2}$:
\begin{equation}
  \label{eq:17}
  \rho_1(f)\,\Tilde\cup\, \rho_1(g) =
  \bigl( d\log f \wedge d\log g,
  -\pi_1(d\log f)\, \log\,\abs{g}
  +\log\,\abs{f}\, \pi_1(d\log g) \bigr)\,.
\end{equation}
The first term is obviously zero when $X$ is a curve. Given $f$
and $g$, invertible on $X$, let us define the imaginary $1$-form:
\begin{equation}
  \label{eq:18}
  r_2(f,g) = \pi_1(d\log f)\, \log\,\abs{g}
  -\log\,\abs{f}\, \pi_1(d\log g)\,.
\end{equation}
On the other hand, the cup product of $f$ and $g$ as elements of
$\delH[1]{X}{\ZZ}{1}$ yields an element $f\cup g$ of
$\delH[2]{X}{\ZZ}{2}$ represented by
\begin{math}
  \bigl( d\log f \wedge d\log g, \log f\; d\log g \bigr)
\end{math}
(if we use the $\cup_0$ product) and this maps via $\rho_2$ to
the element
\begin{equation*}
  \bigl( d\log f \wedge d\log g, -\pi_1(\log f\; d\log g)
  \bigr)\,.
\end{equation*}
This is equal to~\eqref{eq:18} only up to homotopy. Indeed, using
\begin{math}
  \pi_{p+q-1}(a\wedge b) = \pi_p(a)\wedge \pi_{q-1}(b) +
  \pi_{p-1}(a)\wedge \pi_q(b)\,,
\end{math}
we find
\begin{equation}
  \label{eq:19}
  r_2(f,g) = d \bigl( \pi_1(\log f)\, \log\,\abs{g} \bigr)
  -\pi_1 \bigl( \log f\; d\log g \bigr)\,,
\end{equation}
where the first term is just the explicit homotopy as computed
from~\eqref{eq:16}.

Recall that the \emph{tame symbol}
(\cite{del:symbole,bloch:dilog_lie}) $\tame{f}{g}$ associated to
$f$ and $g$ is the line bundle with connection determined (up to
isomorphism) by the class $f\cup g$. A ``Bloch dilogarithm''
(\cite{esn-vie:del}) is (the logarithm of) a horizontal
trivializing section, namely a function $L$ on $U\subset X$
satisfying the equation
\begin{equation*}
  dL = -\log f\; d\log g\,.
\end{equation*}
Thus $L$ realizes the isomorphism $\tame{f}{g}\cong \sho{X}$ over
$U\subset X$.

Thus a Bloch dilogarithm will only locally be available, in
general. However, when $g=1-f$ then $\tame{1-f}{f}$ is
\emph{globally} trivial~\cite{bloch:dilog_lie,esn-vie:del}, i.e.
$\tame{1-f}{f}\cong \sho{X}$. This is the the \emph{Steinberg
  relation} satisfied by the Tame symbol. It can be deduced from
the following universal case. Set $f=z$, $X=\PP^1 \setminus
\{0,1,\infty\}$, consequently $g=1-z$. Then $L$ is identified
with the classical Euler dilogarithm $\li$, namely
\begin{equation*}
  \li (z) = -\int_0^z \log (1-t) \frac{dt}{t}\,,
\end{equation*}
see~\cite{del:symbole} and~\cite{MR94k:19002} for details.

On $\PP^1 \setminus \lbrace 0,1,\infty\rbrace$ the classical
dilogarithm has a single valued parter, denoted $\bwd$,
introduced by Bloch and Wigner:
\begin{equation}
  \label{eq:20}
  \bwd (z) = \arg (1-z) \log\,\abs{z} + \im\li (z)\,.
\end{equation}
$\bwd$ is real-analytic on $\PP^1 \setminus \{0,1,\infty\}$ and
extends continuously to $\PP^1$. That it is single-valued can be
verified as follows. Choose a branch of the logarithm, say the
principal one, to define $\li$ (and $\bwd$) on the cut plane $\CC
\setminus (-\infty,0] \cup [1,\infty)$. Then one shows that the
expression~\eqref{eq:20} is in fact single-valued by analytic
continuation along paths based, say, at $1/2\in \CC$, and winding
around the points $0$ and $1\in\CC$. Explicit computations can be
found in ref.~\cite{MR2001i:11082}.

It is convenient to introduce
\begin{equation*}
  \bwli (z) = \sqrt{-1}\, \bwd (z) \,,
\end{equation*}
so that
\begin{equation*}
  d\bwli = r_2(1-z,z)\,.
\end{equation*}
More generally, if $L$ trivializes $\tame{f}{g}$ over $U$ in the
sense explained above, we can associate a function $\bwli (f,g)$
over $U$ such that
\begin{equation*}
  d \bwli(f,g) = r_2(f,g)
\end{equation*}
via the position
\begin{equation}
  \label{eq:21}
  \bwli (f,g) = \pi_1(\log f)\, \log\,\abs{g} +\im L\,.
\end{equation}

\section{Constructions in hermitian holomorphic Deligne cohomology}
\label{sec:herm-holom-deligne}

In this section we recall the definition of hermitian holomorphic
Deligne cohomology. In ref.~\cite{bry:quillen} Brylinski
introduced certain complexes $C(l)^\bullet$, for a positive
integer $l$, in order to compare the \bei{}-Chern classes of a
holomorphic vector bundle $E$ with the Cheeger-Chern-Simons
classes determined by $(E,\nabla)$, where $\nabla$ is the
\emph{canonical connection,} namely the unique connection
compatible with both the holomorphic and hermitian structures.
The cohomology groups determined by these complexes are aptly
named \emph{Hermitian Holomorphic Deligne cohomology groups.} For
a holomorphic \emph{line} bundle equipped with the canonical
connection, the complex $C(1)^\bullet$ encodes the reduction of
the structure group from $\CC$ to $\TT$ afforded by the hermitian
fiber metric.

In the following we will need to compute Hermitian Holomorphic
Deligne cohomology by means of different---but
quasi-isomorphic---sheaf complexes we denote $\dhh{l}$. These
complexes are tailored to a direct description of a metrized line
bundle in terms of local representatives of the hermitian fiber
metric. Since the two constructions are related by a
quasi-isomorphism, the resulting cohomology groups are the same.

\subsection{Hermitian holomorphic Deligne cohomology}
\label{sec:herm-holom-deligne-1}

In ref.~\cite{bry:quillen}, where Brylinski introduces the
complexes:
\begin{equation}
  \label{eq:24}
  C(l)^\bullet = \cone \bigl(
  \ZZ(l)_X \oplus (F^l\!\sha{X}\cap \sigma^{2l}\she{X}(l))
  \longrightarrow \she{X}(l)
  \bigr)[-1]\,,
\end{equation}
where $\sigma^p$ denotes the (sharp) truncation in degree $p$,
namely for a complex $\sheaf{F}^\bullet$, $\sigma^p\sheaf{F}^k$
is zero for $k<p$ and equal to $\sheaf{F}^k$ for $k\geq p$.
In~\eqref{eq:24} we take the cone of the difference between the
two inclusions. We have the following:
\begin{definition}[\cite{bry:quillen}]
  \label{def:1}
  The hermitian holomorphic Deligne cohomology groups are the
  hypercohomology groups of the complexes~\eqref{eq:24}:
  \begin{equation}
    \label{eq:23}
    \dhhH[p]{X}{l} \eqdef \HHH^{\,p}(X,C(l)^\bullet)\,.
  \end{equation}
\end{definition}
The complexes~\eqref{eq:24} are expressed as cones, and therefore
admit a (standard) \bei\ product~\eqref{eq:2}. (The wedge product
induces cup products on both $\she{X}(l)$ and $F^l\!\sha{X}\cap
\sigma^{2l}\she{X}(l)$.)\footnote{It appears signs should be
  adjusted in the product formula quoted in
  ref.~\cite{bry:quillen}, and that using~\eqref{eq:2} is more
  appropriate.} It follows there is a graded commutative product
on cohomology:
\begin{equation}
  \label{eq:27}
  \dhhH[i]{X}{l} \otimes \dhhH[j]{X}{k}
  \overset{\cup}{\longrightarrow}
  \dhhH[i+j]{X}{l+k}\,.
\end{equation}
Also, from the standard cone exact sequences we get
(cf.\ ref.~\cite{bry:quillen}):
\begin{equation}
  \label{eq:84}
  \dotsm \longrightarrow H^{2l-1}(X,\RR (l))
  \longrightarrow \dhhH[2l]{X}{l} \longrightarrow
  H^{2l}(X,\ZZ (l)) \oplus A^{(l,l)}(X)_{\RR(l)}
  \longrightarrow H^{2l}(X,\RR (l))
  \longrightarrow \dotsm
\end{equation}
where $A^{(l,l)}(X)_{\RR(l)}$ denotes the space of smooth
$\RR(l)$-valued global $(l,l)$-forms on $X$. Thus we see
hermitian holormorphic classes are $\ZZ(l)$-valued classes
represented by (necessarily) closed $\RR(l)$-valued $2l$-forms of
pure type $(l,l)$.  For a line bundle this corresponds to a
structure group reduction from $\sho{X}$ to $\she[0]{X}$, namely
$\TT$-valued sections, at the same time controlling the Hodge
type of the resulting class, cf.\ ref.~\cite{bry:quillen}

Later (cf. sec.~\ref{sec:line-bun}) we will want to work with the
hermitian structure on a holomorphic line bundle, together with
the imaginary $(1,1)$-form built from the canonical connection,
directly in the holomorphic frame. To carry this out in general
for $(l,l)$ classes, we introduce the complex
\begin{equation}
  \label{eq:22}
  \dhh{l} = \cone \bigl(
  \deligne{\ZZ}{l}\oplus (F^l\!\sha{X}\cap \sigma^{2l}\she{X}(l))
  \longrightarrow \deltilde{\RR}{l}
  \bigr)[-1]\,.
\end{equation}
The map
\begin{math}
  \deligne{\ZZ}{l} \to \deltilde{\RR}{l}
\end{math}
is the composite of the obvious map
\begin{math}
  \deligne{\ZZ}{l} \to \deligne{\RR}{l}
\end{math}
with the quasi-isomorphism $\rho_{\,l}$ defined by~\eqref{eq:12}. We
will simply denote it by $\rho_l$ in the sequel, suppressing the first
morphism in the notation. The map
\begin{math}
  (F^l\!\sha{X}\cap \sigma^{2l}\she{X}(l)) \to \deltilde{\RR}{l}
  = \cone \bigl( F^l\!\sha{X} \to \she{X}(l-1) \bigr)[-1]
\end{math}
is induced by the inclusion of $(F^l\!\sha{X}\cap
\sigma^{2l}\she{X}(l))$ into $F^l\!\sha{X}$. In~\eqref{eq:22} we
take the cone of the difference between these two maps.  The
complex~\eqref{eq:22} offers another model for Hermitian
holomorphic Deligne cohomology. Indeed we have:
\begin{lemma}
  \label{lemma:qi-dhh-bry}
  The complexes $C(l)^\bullet$ and $\dhh{l}$ are
  quasi-isomorphic.
\end{lemma}
\begin{proof}
  By elementary manipulation of cones
  \begin{equation*}
    \dhh{l}= \cone \Bigl(
    F^l\!\sha{X}\cap \sigma^{2l}\she{X}(l) \to
    \cone \bigl( \deligne{\ZZ}{l} \to
    \deltilde{\RR}{l} \bigr) \Bigr)[-1]
  \end{equation*}
  and clearly:
  \begin{math}
    \cone \bigl( \deligne{\ZZ}{l} \to \deltilde{\RR}{l} \bigr)
    \qi \cone \bigl( \ZZ(l)_X \to \RR(l)_X \bigr) \qi \cone
    \bigl( \ZZ(l)_X \to \she{X}(l) \bigr)\,,
  \end{math}
  where all arrows are quasi-isomorphisms. Thus
  \begin{align*}
    \dhh{l} &\lqi \cone \Bigl( F^l\!\sha{X}\cap
    \sigma^{2l}\she{X}(l) \to
    \cone \bigl( \ZZ(l)_X \to \she{X}(l) \bigr) \Bigr)[-1] \\
    &= \cone \bigl( \ZZ(l)_X \oplus (F^l\!\sha{X}\cap
    \sigma^{2l}\she{X}(l)) \to \she{X}(l) \bigr)[-1] \\ &\equiv
    C(l)\,,
  \end{align*}
  as wanted.
\end{proof}
It follows from Lemma~\ref{lemma:qi-dhh-bry} that
\begin{math}
  \HHH^p(X,\dhh{l}) \cong \dhhH[p]{X}{l}\,,
\end{math}
so we can use either complex to compute the Hermitian-holomorphic
Deligne cohomology groups. 

Again from the cone exact sequence applied to~\eqref{eq:22}, we
see the groups $\dhhH{X}{l}$ also satisfy the exact sequence
\begin{equation}
  \label{eq:25}
  \dotsm \longrightarrow \delH[2l-1]{X}{\RR}{l}
  \longrightarrow
  \dhhH[2l]{X}{l} \longrightarrow
  \delH[{2l}]{X}{\ZZ}{l} \oplus A^{(l,l)}(X)_{\RR(l)}
  \longrightarrow \delH[{2l}]{X}{\RR}{l}
  \longrightarrow \dotsm
\end{equation}
which we can rewrite as
\begin{equation*}
  \dotsm \longrightarrow \delH[2l-1]{X}{\RR}{l}
  \longrightarrow \dhhH[2l]{X}{l} 
  \longrightarrow \delH[{2l}]{X}{\ZZ}{l} \oplus
  A^{(l,l)}(X)_{\RR(l)} \longrightarrow H^{2l}(X,\RR(l))\cap
  H^{l,l}(X) \longrightarrow \dotsm
\end{equation*}
Thus the elements of $\dhhH[2l]{X}{l}$ map onto those
$(l,l)$-forms representing the Hodge classes corresponding to
$\delH[2l]{X}{\ZZ}{l}$.

The complexes $\deligne{\ZZ}{l}$, $\deltilde{\RR}{l}$, and
$F^l\!\sha{X}\cap \sigma^{2l}\she{X}(l)$ appearing in the
cone~\eqref{eq:22}, all have cup products. It follows that we
have the \bei{} family of products
\begin{equation}
  \label{eq:26}
  \dhh{l} \otimes \dhh{k} \overset{\cup_\alpha}{\longrightarrow} 
  \dhh{l+k} \,.
\end{equation}
\begin{remark}
  The product $\cup_\alpha$ in eq.\eqref{eq:26} should be
  intended in the modified version provided by~\eqref{eq:5},
  since the complexes $\deligne{\ZZ}{l}$ and $\deltilde{\RR}{l}$
  have product structures that are compatible with the map
  $\rho_l$ only up to homotopy (given by formula~\eqref{eq:16}).
  Moreover, these complexes have product structures that are
  themselves graded commutative up to homotopy: that the
  product~\eqref{eq:26} is graded commutative up to homotopy
  follows from Proposition~\ref{prop:4} in the appendix.
\end{remark}

\subsection{Hermitian holomorphic line bundles}
\label{sec:line-bun}

A hermitian holomorphic line bundle or, equivalently, a metrized
line bundle, cf. \cite{lang:arakelov}, is a holomorphic line
bundle $L$ over $X$ together with a hermitian fiber metric $\rho
: L \rightarrow \RR_{\geq 0}$. As a rule, we will not distinguish
$L$ and its sheaf of holomorphic sections. We will also use the
alternate notation $\abs{s}_\rho$ to denote the length of a local
section $s$ of $L$ with respect to $\rho$. Metrized line bundles
can be tensor multiplied and an inverse is defined, see. op. cit.
An isomorphism of metrized line bundles $(L,\rho)$ and
$(L^\prime,\rho^\prime)$ is defined in the obvious way, namely it
is a map $\phi : L \to L^\prime$ such that $\abs{s}_\rho =
\abs{\phi(s)}_{\rho^\prime}$ for some local section $s$ of $L$.
We denote by $\widehat{\pic{X}}$ the group of isomorphism classes
of metrized line bundles.

If $L$ is trivialized over a \cech\ cover $\cover{U}_X =
\{U_i\}_{i\in I}$ by sections $s_i$, then as usual we obtain the
cocycle of transition functions $g_{ij}\in \sho{X}^\times
(U_i\cap U_j)$ by writing $s_j = s_i\, g_{ij}$.  Then if
$(L,\rho)$ is a metrized line bundle, we can define the positive
function $\rho_i = \abs{s_i}^2_\rho$, namely the local
representative of the hermitian structure with respect to the
given trivialization. It follows that the various local
representatives satisfy
\begin{equation}
  \label{eq:29}
  \rho_j = \rho_i \, \abs{g_{ij}}^2\,.
\end{equation}
Let us work out the local version of the isomorphism introduced
above. Let $s'_i$ be a local section of $L'$ over $U_i$.
Introduce analogous (primed) quantities for $L'$ as we just did
for $L$. Given the isomorphism $\phi: L\to L'$ we have $\phi(s_i)
= s'_i\, f_i$, for some $f_i\in \sho{X}^\times(U_i)$. Then we
find $f_i\, g_{ij} = g'_{ij}\, f_j$ and $\rho_i = \rho'_i\,
\abs{f_i}^2$.

Still working with respect to the chosen cover $\cover{U}_X$, a
\emph{connection} compatible with the holomorphic structure is
the datum of a collection of $(1,0)$-forms $\xi_i \in
\sha[{1,0}]{X}(U_i)$ satisfying
\begin{equation}
  \label{eq:30}
  \xi_j - \xi_i = d\log g_{ij}\,.
\end{equation}
Note for future reference that $\sha[{1,0}]{X} =
F^1\!\sha[1]{X}$. The connection is compatible with the hermitian
metric if
\begin{equation}
  \label{eq:31}
  \pi_0(\xi_i) = \onehalf d\log \rho_i\,.
\end{equation}
Using $d=\del +\delb$ and decomposition with respect to
$(p,q)$-types, we find the familiar relation
\begin{equation}
  \label{eq:32}
  \xi_i = \del \log \rho_i
\end{equation}
for the \emph{unique} connection compatible with both the complex
and hermitian structures \cite{gh:alg_geom}.

The global $2$-form
\begin{equation}
  \label{eq:33}
  c_1(\rho) = \delb\del\log \rho_i
\end{equation}
represents the first Chern class of $L$ in $H^2(X,\RR(1))$.
Actually, the class of $c_1(\rho)$ is a pure Hodge class in
$H^{1,1}(X)$ and, according to the examples, it is the image of
the first Chern class of $L$ under the map $\delH[2]{X}{\ZZ}{1}
\to H^2_\mathcal{D}(X,\RR(1))$ induced by $\ZZ(1) \to \RR(1)$.
Observe that $c_1(\rho) = c_1(\rho')$ under the isomorphism
considered above.

The following proposition can be found, for example, in
refs.~\cite{brymcl:deg4_II,bry:quillen}. (It apparently was
observed first by Deligne, cf.~\cite{esn:char}.) The proof is
based on writing out an explicit cocycle in a \cech\ resolution
of $C(1)^\bullet$ or $\dhh{1}$. This calculation will be needed
later on with the complex $\dhh{1}$, so we provide a proof here.
\begin{proposition}
  \label{prop:1}
  \begin{equation*}
      \widehat{\pic{X}} \cong \dhhH[2]{X}{1}
  \end{equation*}
\end{proposition}
\begin{proof}
  Recall that $\dhh{1}$ is the cone of the map $\rho_1 -\jmath$,
  where $\rho_1: \deligne{\ZZ}{1} \to \deltilde{\RR}{1}$, and
  \begin{math}
    \jmath : F^1\!\sha{X}\cap \sigma^{2}\she{X}(1) \to
    \deltilde{\RR}{1}\,.
  \end{math}
  By unraveling the structure of all the cones involved we have:
  \begin{equation*}
    \begin{CD}
      \ZZ(1)_X @>>> \shomega[1]{X}\oplus\sho{X} @>>>
      \shomega[2]{X}\oplus \shomega[1]{X} @>>>
      \shomega[3]{X}\oplus \shomega[2]{X} @>>> \dotsm \\
      & & @VV{\imath\oplus\pi_0}V @VV{\imath\oplus\pi_0}V
      @VV{\imath\oplus\pi_0}V \\
      & & F^1\!\sha[1]{X} \oplus \she[0]{X} @>>> F^1\!\sha[2]{X}
      \oplus \she[1]{X} @>>>
      F^1\!\sha[3]{X} \oplus \she[2]{X} @>>> \dotsm \\
      & & & & @AA{\jmath\oplus 0}A @AA{\jmath\oplus 0}A \\
      & & & & F^1\!\sha[2]{X}\cap \she[2]{X}(1) @>>>
      F^1\!\sha[3]{X}\cap \she[3]{X}(1) @>>> \dotsm
    \end{CD}
  \end{equation*}
  With respect to this diagram, an element of total degree $2$
  can be written in the form:
  \begin{equation}
    \label{eq:37}
    \setlength{\extrarowheight}{4pt}
    \begin{array}{c|c|c}
      c_{ijk} & -\;d \log\,g_{ij}\oplus \log\,g_{ij} & 0 \\
      \hline
      0 & \xi_i \oplus \sigma_i & \mathsf{X} \\
      \hline
      0 & 0 & \eta_i
    \end{array}
  \end{equation}
  for $\xi_i\in F^1\!\sha[1]{X}(U_i)$, $\sigma_i \in
  \she[0]{X}(U_i)$ and $\eta_i\in (F^1\!\sha[2]{X}\cap
  \she[2]{X}(1))(U_i)$. To make sense out of~\eqref{eq:37}, note
  that each entry is an element of the object in the
  corresponding position in the left $3\times 3$ part of the
  previous diagram. Then since the total degree is $2$, the
  degree of each element \emph{in the complex} $\dhh{1}$ is $2$
  minus the \cech\ degree as found in~\eqref{eq:37}. The top line
  is the class of $L$ in $\deligne[2]{\ZZ}{1}$. Finally, the
  entry marked $\mathsf{X}$ means there is no applicable
  element---it would have degree $3$.
  
  A totally routine calculation shows that~\eqref{eq:37} is a
  degree $2$ cocycle if and only if
  relations~\eqref{eq:30},~\eqref{eq:31} and $\eta_i = d\xi_i =
  2\,\delb\del \sigma_i$ are satisfied (with $\sigma_i = \onehalf
  \log \rho_i$). Thus $\eta_i = c_1(\rho)\vert_{U_i}$ and all the
  relations defining a metrized line bundle with its canonical
  connection are satisfied.  The verification that adding an
  appropriate coboundary to the cocycle leads to an isomorphic
  metrized bundle in the sense explained above is also routine.
  Finally, the correspondence between component-wise addition of
  cocycles modulo coboundaries and the group operation in
  $\widehat{\pic{X}}$ is again a direct verification.
\end{proof}
From the proof of Proposition~\ref{prop:1} we have that the
cocycle representing $(L,\rho)$ in the \cech\ resolution of
$\dhh{1}$ has the expression:
\begin{equation}
  \label{eq:38}
  \setlength{\extrarowheight}{4pt}
  \begin{array}{c|c|c}
    c_{ijk} & -\;d \log\,g_{ij}\oplus \log\,g_{ij} & 0 \\
    \hline
    0 & \del \log\rho_i \oplus \onehalf \log \rho_i & \mathsf{X} \\
    \hline
    0 & 0 & \delb\del\,\log \rho_i 
  \end{array}
\end{equation}
This cocycle is expressed purely in terms of holomorphic and
metric data (the local expression for the fiber metric $\rho$),
and it also explicitly encodes the canonical connection and its
curvature.

By way of comparison, a degree two cocycle in the \cech\ 
resolution of $C(1)^\bullet$ with respect to the same cover would
be given by a triple
\begin{equation*}
  \bigl( \eta_i\oplus \sigma_i, \log u_{ij}, c_{ijk} \bigr)\,,
\end{equation*}
where $\eta_i\in \sha[(1,1)]{X}(U_i)_{\RR(1)}$, $\sigma_i \in
\she[1]{X}(1)(U_i)$, $\log u_{ij}\in \she[0]{X}(1)(U_{ij})$ and
$c_{ijk}\in \ZZ(1)$, satisfying the following relations: aside
from the standard cocycle condition for the $\TT$-valued
functions $u_{ij}$, we must have
\begin{equation}
  \label{eq:86}
  \begin{gathered}
    d\eta_i = 0\,,\\
    d\sigma_i=\eta_i\,,\\
    \sigma_j -\sigma_i = d \log u_{ij}\,.
  \end{gathered}
\end{equation}
Given a holomorphic trivialization of $(L,\rho)$ as above, the
position
\begin{equation}
  \label{eq:87}
  \eta_i = \delb\del\log \rho_i\,, \quad
  \sigma_i = \onehalf d^c \log \rho_i\,, \quad
  u_{ij} = g_{ij}/\abs{g_{ij}}
\end{equation}
satisfies relations~\eqref{eq:86}. From here, the canonical
connection is recovered. It is also easily seen that another
solution of~\eqref{eq:86} is equivalent to the one above. The
following is evident:
\begin{lemma}
  \label{lem:1}
  Relations~\eqref{eq:87} provide an explicit quasi-isomorphism
  between the \cech\ resolutions of $\dhh{1}$ and $C(1)^\bullet$.
\end{lemma}

On the
other hand, in~\eqref{eq:38} the relations defining the canonical
connection are already enforced at the level of the complex.

Finally, if $[L,\rho]$ and $[L',\rho']$ are the classes
corresponding to the metrized bundles $(L,\rho)$ and
$(L',\rho')$, then we write
\begin{math}
  [ L\otimes L', \rho\rho'] = [L,\rho] + [L',\rho']\,.
\end{math}

\subsection{Cup product of hermitian holomorphic line bundles}
\label{sec:Cup-prod-herm}

If $L$ and $L'$ are two line bundles on $X$, their cup product in
Deligne cohomology would produce a class in $\delH[4]{X}{\ZZ}{2}$
denoted $\tame{L}{L'}$. Calculations were carried out
in~\cite{brymcl:deg4_II} where a geometric interpretation of
$\tame{L}{L'}$ as a 2-gerbe and its connection with the
determinant of cohomology (when $X$ is the total space of a
family of Riemann surfaces) were put forward. The structure
ensuing from the generalization to line bundles equipped with a
hermitian metric was further analyzed in ref.~\cite{bry:quillen}
by means of hermitian holomorphic Deligne cohomology.

The cup product of two metrized line bundles $(L,\rho)$ and
$(L',\rho')$ in hermitian holomorphic Deligne cohomology produces
a class in $\dhhH[4]{X}{2}$. Despite the more involved definition
of $\dhh{l}$ as opposed to that of $C(l)^\bullet$, it will
nonetheless be more advantageous from the perspective of
sect.~\ref{sec:conf-metr-herm} to use the former to calculate the
desired class. Thus let us explicitly compute a representative of
this class in terms of the expression~\eqref{eq:38} for the class
of a metrized line bundle and the modified \bei\ 
product~\eqref{eq:5} for
\begin{equation*}
  \dhh{1} \otimes \dhh{1} \xrightarrow{\cup_{\alpha}}
  \dhh{2}  \,. 
\end{equation*}
A diagram for the complex $\dhh{2}$ analogous to that for
$\dhh{1}$ displayed in the proof of Proposition~\ref{prop:1} is:
\begin{equation}
  \label{eq:39}
  \begin{CD}
    \ZZ(2)_X @>>> \sho{X} @>>> \shomega[2]{X}\oplus
    \shomega[1]{X} @>>> \shomega[3]{X}\oplus \shomega[2]{X} @>>>
    \shomega[4]{X}\oplus \shomega[3]{X} @>>>
    \dotsm \\
    & & @VV{\pi_1}V @VV{\imath\oplus\pi_1}V
    @VV{\imath\oplus\pi_1}V @VV{\imath\oplus\pi_1}V \\
    & & \she[0]{X}(1) @>>> F^2\!\sha[2]{X} \oplus \she[1]{X}(1)
    @>>> F^2\!\sha[3]{X} \oplus \she[2]{X}(1) @>>>
    F^2\!\sha[4]{X} \oplus \she[3]{X}(1) @>>>
    \dotsm \\
    & & & & & & & & @AA{\jmath\oplus 0}A \\
    & & & & & & & & F^1\!\sha[4]{X}\cap \she[4]{X}(2) @>>> \dotsm
  \end{CD}
\end{equation}
The cup product of two elements of the form~\eqref{eq:38} in the
\cech\ resolution of $\dhh{1}$ would result in an analogous
$5\times 3$ table. The degrees (in the cone) start from $0$ in
the leftmost entry in the top and bottom rows, and from $1$ in
the center one. The overall signs are determined by those in the
\bei\ product~\eqref{eq:5} plus those arising from the \cech\ 
resolution as explained in sect.~\ref{sec:notation}. Actually, it
is visually less cumbersome to display the resulting expression
as in~\eqref{tab:big-cup} below, where the corresponding
bidegrees are explicitly indicated: the first degree is the
overall degree in the cone and the second is the \cech\ degree.
(We have explicitly written only the nonzero terms.)
\begin{small}
  \begin{equation}
    \renewcommand{\onehalf}{\tfrac{1}{2}}
    \renewcommand{\arraystretch}{1.3}
    \begin{array}{c|c}
      (0,4) & c_{ijk}c'_{klm} \\ \hline
      (1,3) & -c_{ijk} \log \, g'_{kl} \\ \hline
      (2,2) & \bigl( -d \log \, g_{ij} \wedge d \log\, g'_{jk} \oplus
      -\log \, g_{ij} \, d\log\,g'_{jk} \bigr) \bigoplus
      -\pi_1\log\,g_{ij}\, \log\,\abs{g'_{jk}} \\ \hline
      (3,1) & 
      \begin{aligned}
        \bigl( -(1-\alpha) \, d\log g_{ij} \wedge \del \log
        \rho'_j +\alpha \, \del \log \rho_i \wedge d \log g'_{ij}
        \bigr)  
        \oplus 
        \Bigl\{ (1-\alpha) \bigl( d^c\log\,\abs{g_{ij}}\;
        \onehalf \log\rho'_j &-\log\,\abs{g_{ij}}\; \onehalf d^c
        \log\,\rho'_j
        \bigr) \\
         +\alpha \bigl( -\onehalf d^c
        \log\rho_i\;\log\,\abs{g'_{ij}} &+\onehalf \log\rho_i\;
        d^c\log\,\abs{g'_{ij}} \bigr) \Bigr\}
      \end{aligned} \\ \hline
      (4,0) & 
      \begin{aligned}
        \delb\del \log\rho_i \wedge
        \delb\del \log \rho'_i
        \bigoplus
        \bigl( & \alpha\, \delb\del\log\rho_i \wedge \del \log
        \rho'_i +(1-\alpha)\, \del\log \rho_i \wedge
        \delb\del\log\rho'_i
        \bigr) \\
        & \oplus \bigl( \alpha\, \delb\del \log\rho_i\;
        \onehalf\log\rho'_i +(1-\alpha)\, \onehalf \log\rho_i \;
        \delb\del\log\rho'_i \bigr)
      \end{aligned}      
    \end{array}
    \label{tab:big-cup}
  \end{equation}
\end{small}

\subsection{Integration along the fiber}
\label{sec:integr-along-fiber}

To conclude this introduction, let us quickly mention that
Hermitian holomorphic Deligne cohomology has an integration along
the fiber map. Namely, if $\pi\colon X\rightarrow S$ is a proper
submersion of complex manifolds, it follows from
\cite{MR91j:14017} that Deligne cohomology has an integration
along the fiber $\int_\pi$, hence by
\cite{brymcl:deg4_II,bry:quillen} there is a map
\begin{equation}
  \label{eq:28}
  \dhhH[i]{X}{l} \longrightarrow \dhhH[i-2d]{S}{l-d}
\end{equation}
where $d$ is the complex dimension of the fiber, and commutative
diagrams analogous to \cite[Theorem 5.1]{bry:quillen}.

We will be interested in the case of complex relative dimension
$1$, namely $\pi \colon X\rightarrow S$ is a holomorphic
fibration with compact connected Riemann surfaces as fibers. From
\eqref{eq:28} in degree $4$ we have the map:
\begin{math}
  \dhhH[4]{X}{2} \rightarrow \dhhH[2]{S}{1}
\end{math}
into the group of complex hermitian line bundles on $S$ (cf.\ 
Proposition \ref{prop:1}).

We refer to refs.\ \cite{brymcl:deg4_II,bry:quillen} for
a complete treatment of the map~\eqref{eq:28}, in particular for
the ``trace'' map
\begin{equation*}
  \RRR^\bullet\pi_* (\dhh{l}) \longrightarrow
  \dhh[\bullet-2d]{l-d}\,,
\end{equation*}
which induces~\eqref{eq:28} at the level of cohomology. (Clearly
any quasi-isomorphic model for Hermitian-holomorphic cohomology
will do.) We will limit ourselves to observe the following: if
$\cover{U}_S$ is a good cover of $S$, and $\cover{U}_X$ a good
cover of $X$ which refines $\pi^{-1}\cover{U}_S$, the $i$-{th}
direct image $\RRR^i\pi_* \sheaf{F}^\bullet$ of any complex
$\sheaf{F}^\bullet$ on $X$ is the sheaf on $S$ associated to the
presheaf
\begin{equation*}
  V \longmapsto H^i \bigl(\Tot
  \Check{C}^\bullet ( \cover{U}_X\cap \pi^{-1}(V),\sheaf{F}^\bullet)
  \bigr)\,. 
\end{equation*}
Then unraveling the cone structure of ${\dhh{l}}_X$ (or of
$C(l)^\bullet_X$) reduces to computing the direct images of
$\ZZ(l)_X$ and $\sha{X}$, the de~Rham complex. The trace map is
then obtained by capping a total cocycle, say in
\begin{equation}
  \label{eq:89}
  \bigoplus_{p+q=i}
  \Check{C}^q ( \cover{U}_X\cap \pi^{-1}(V),\sha[p]{X})\,,
\end{equation}
with a representative of the fundamental class $[M]$ of the
smooth model $M$ of the fiber of $\pi^{-1}(V)\cong V\times M$ (as
smooth manifolds). Explicit representatives for $[M]$ for a
triangulation subordinated to the nerve of $\cover{U}_X\cap
\pi^{-1}(V)$ are computed in ref.\ \cite{aldtak2000}.  If
$\omega^p_{i_0,\dots,i_q}$ are the components of a cocycle
$\omega^i_V$ in~\eqref{eq:89}, then
\begin{equation*}
  \omega^i_V \cap [M]
  =\sum \omega^p_{i_0,\dots,i_q} \cap \Delta^{2d-q}_{i_0,\dots,i_q}
\end{equation*}
where the $\Delta^{2d-q}_{i_0,\dots,i_q}$ are \emph{signed}
generators of the nerve of dimension $2d -q$. The exact
combinatorics can be found in ref.\ \cite{aldtak2000}, so we will
not further pursue the matter here.

Let us conclude by noticing that for $\pi\colon X\to S$ the cup
product of $(L,\rho)$ and $(L',\rho')$ induces a metrized line
bundle on $S$. {F}rom the preceding discussion, the connection
component, for example, will be obtained by following the (two)
components $F^2\!\sha{X}\cap \she[\bullet \geq 4]{X}(2)$
in~\eqref{tab:big-cup}. The one in degree $(4,0)$ will be capped
with generators $\Delta^{2}$ of dimension $2$ (hence integrated
over $2$-simplices), those in degree $(3,1)$ with generators
$\Delta^{1}$ of dimension $1$.

From \cite[Theorem 6.1]{brymcl:deg4_II} there is a
corresponding map
\begin{math}
  A^{(2,2)}(X)_{\RR(2)} \rightarrow A^{(1,1)}(S)_{\RR(1)}\,.
\end{math}
{F}rom the explicit cocycles we have computed the corresponding
representative in $H^2(S,\RR(1))\cap H^{1,1}(S)$ of the first
Chern class of the resulting metrized line bundle is
\begin{equation*}
  \int_\pi c_1(\rho) \wedge c_1(\rho')\,,
\end{equation*}
as in \cite[Proposition 6.6.1]{MR89b:32038}.

\section{Conformal metrics and the Liouville functional}
\label{sec:Conf-metr-Liouv}

Let $X$ be a compact Riemann surface of genus $g\geq 2.$ For simplicity
we can assume $X$ to be connected, although this is not necessary, and
in fact this assumption will be dropped when dealing with Kleinian
groups.

Let $\cm{X}$ be the space of conformal metrics on $X$. Locally on $X$,
if $z$ is a local analytic coordinate defined on an open set $U$, any
metric $ds^2$ can be represented as
\begin{equation*}
  ds^2 = \rho \, \abs{dz}^2
\end{equation*}
for a positive function $\rho : U \to \RR_{>0}$. According to
sec.~\ref{sec:line-bun}, a conformal metric corresponds to
considering the metrized line bundle $(T_X,\rho)$, where $T_X$ is
the holomorphic tangent bundle of $X$. With respect to a cover
$\cover{U}_X$ of $X$, the conformal factors $\rho_i$ and
$\rho_j$, associated to $U_i$ and $U_j$ respectively, satisfy the
relation
\begin{equation}
  \label{eq:41}
  \rho_j = \rho_i\, \abs{z'_{ij}}^2\,,
\end{equation}
where $z'_{ij} = dz_i/dz_j$ and $z_i$ is a local coordinate
defined over $U_i$.

It follows from the uniformization theorem that there exists a
unique conformal metric of scalar curvature equal to $-1$, the
Poincar\'e metric. Locally on $U_i\subset X$, the condition for
the metric to have curvature $-1$ is equivalent to the nonlinear
PDE
\begin{equation}
  \label{eq:42}
  \frac{\del^2}{\del z_i\del\Bar z_i} \phi_i = \onehalf
  \exp {\phi_i} \,,
\end{equation}
known as the Liouville equation, for the smooth function $\phi_i
= \log \rho_i$. Observe that equation~\eqref{eq:42} can be
written in the form
\begin{equation*}
  c_1(\rho) = \sqrt{-1}\,\omega_\rho\,,
\end{equation*}
where we have used the K\"ahler form associated to the metric:
\begin{equation*}
  \omega_\rho\vert_{U_i} = \ihalf \rho_i \,
  dz_i\wedge d\Bar z_i\,.  
\end{equation*}
This representation makes it apparent that the Liouville equation
is independent of the choice of the coordinate system. On the
other hand, a direct verification of this fact is immediate using
$\phi_j -\phi_i = \log \abs{z'_{ij}}^2$.

\subsection{Variational problem for conformal metrics}
\label{sec:vari-probl-conf}

It is well known that the Liouville equation has a local
variational principle in the following sense. Let $D$ be a region
in the complex plane. Then equation~\eqref{eq:42} is the
Euler-Lagrange equation for the variational problem defined by
the action functional
\begin{equation}
  \label{eq:43}
  S[\phi] = \ihalf \int_D \bigl( \del \phi \wedge \delb \phi
  + e^\phi\,dz\wedge d\Bar z \,\bigr) \,,
\end{equation}
defined on smooth functions $\phi :D \to \RR$, with the condition
that variations
\begin{math}
  \left.\frac{d}{d\alpha} (\phi_\alpha) \right\vert_{\alpha=0}
\end{math}
of $\phi$ be zero on $\del D$. However, it is easily seen that
the functional~\eqref{eq:43} cannot be defined globally on $X$,
since, as a consequence of~\eqref{eq:41}, the first term under
the integral sign would not yield a well-defined $2$-form on $X$.
(The second one would of course present no problems, it would
just give the area $A_X(\rho)$ of $X$ with respect to the given
metric $\rho\,\abs{dz}^2 = e^\phi\,\abs{dz}^2$.)
Accordingly, it is convenient to write the integrand in~\eqref{eq:43} as 
\begin{equation*}
  \sqrt{-1}\,\omega[\phi] + \omega_\rho
\end{equation*}
where we have defined the $2$-form\footnote{Note that
  equation~\eqref{eq:44} defines an imaginary form. The reason for this
  choice will be apparent later.}
\begin{equation}
  \label{eq:44}
  \omega^{0}[\phi] = \onehalf \del \phi \wedge \delb \phi\,,
\end{equation}
and restrict our considerations
to the first term of~\eqref{eq:43} which we denote
\begin{equation}
  \label{eq:45}
  \Check{S}[\phi] = \sqrt{-1}\,\int_D \omega^{0}[\phi]\,.
\end{equation}

There is by now an established procedure on how to address the problem
caused by the fact that~\eqref{eq:44} is not globally defined.  In
general terms, given the choice of a conformal metric $ds^2\in\cm{X}$
and a cover $\cover{U}_X$, one works with the full \cech-de~Rham complex
$\Check{C}^\bullet(\cover{U}_X,\she{X}(1))$ with respect to
$\cover{U}_X$, rather than with just differential forms. The
$2$-form~\eqref{eq:44} is then completed to a total degree
$2$-cocycle---to be denoted $\Omega[\phi]$. This results in a class in
$H^2(X,\RR(1))$ after taking cohomology. (Whether or not there also is a
de~Rham type theorem will depend on the acyclicity properties of
$\cover{U}_X$.)

This scheme has been previously carried out not quite for covers
of $X$ by open sets, but actually for different choices of planar
coverings. For the covering associated to a Schottky
uniformization of $X$ a generalization of eq.~\eqref{eq:43} was
written in ref.\ \cite{zogtak1987-2}. More recently, a detailed
calculation of the cocycle for the general case of a covering
associated to a Kleinian uniformization was carried out in ref.\
\cite{math.CV/0204318} by exploiting the homological methods
developed in ref.\ \cite{aldtak1997}. Note, however, that from the
point of view of ref.\ \cite{aldtak2000} these planar coverings
are ``\'etale'' coordinates on $X$, so the group cohomology
constructions required to work with the various kinds of
uniformization coverings just follow from specializing the \cech\ 
formalism to the coverings at hand.

Finally, the integration in eqs.~\eqref{eq:43} or~\eqref{eq:45} should
be replaced by the evaluation of $\Omega[\phi]$ against an appropriate
representative $\Sigma$ of the fundamental class of $X$. The
``appropriate'' form for both $\Sigma$ and the evaluation will be
dictated by the chosen cover $\cover{U}_X$ \emph{and} the cohomology
theory being used. Typically $\Sigma$ will be a cycle in a double
complex of $\cover{U}_X$-small simplices, where the differentials are
the singular one and the one determined by the face maps induced by the
cover.  Thus, in the case of a \cech\ cover, it will be the complex
determined by a triangulation of $X$ subordinated to the open
$\cover{U}_X$, and in the same way, for a planar cover the singular
complex of the planar domain $\cover{U}_X$ tensored with an appropriate
bar resolution of the group of deck transformations. These issues have
been discussed at length in refs.\ 
\cite{aldtak1997,aldtak2000,math.CV/0204318}, so we will not repeat the
discussion here. Whenever we have a cocycle extending~\eqref{eq:44} and
a cycle $\Sigma$ representing $X$ we state
\begin{definition}
  \label{def:2}
  The Liouville functional (without the area term) is given by
  the evaluation:
  \begin{equation}
    \label{eq:46}
    \Check{S}[\phi] = -\frac{1}{\tate}
    \dual{\Omega[\phi]}{\Sigma}\,. 
  \end{equation}
  For the complete functional we add the area term
  \begin{equation}
    \label{eq:47}
    S[\phi] = \Check{S}[\phi] +\frac{1}{2\pi}
    \int_X\omega_\rho\,. 
  \end{equation}
\end{definition}
\begin{remark}
  $\dual{\Omega[\phi]}{\Sigma}\in \RR(1)$, and $\Check{S}[\phi]$ (or
  $S[\phi]$) is real. Division by $\tate$ is conventional, but note that
  $\RR(1) \qi \RR_+$ via $\exp(\cdot / \tate )$. In the sequel it will
  be more convenient to work with imaginary classes. (See also
  section~\ref{sec:conf-metr-herm}.)
\end{remark}
In next two subsections we examine these constructions in some
detail. For definiteness, we initially make use of an ordinary
open cover. First, we recall the direct construction of a cocycle
generalizing~\eqref{eq:44}, and we show that this way
eq.~\eqref{eq:42} indeed is the resulting extremum condition.
Then we emphasize the role played by Deligne cohomology and the
tame symbol. These aspects will become important when introducing
a Kleinian uniformization later in the paper, when we discuss
connections with the dilogarithm function.  In a later section we
shall tackle the question of its geometrical significance by
making full use of the hermitian-holomorphic version of Deligne
cohomology presented in section~\ref{sec:herm-holom-deligne}, and
we show that the cocycle constructed following refs.
\cite{zogtak1987-2,math.CV/0204318} corresponds to the square of
$(T_X,\rho)$ in hermitian holomorphic Deligne cohomology.

\subsection{Direct construction of the Liouville cocycle}
\label{sec:Direct-constr-Liouv}

\subsubsection{Initial setup}
\label{sec:Initial-setup}

Let $X$ a compact Riemann surface of genus greater than $2$. We
shall not include the area term in our explicit calculations,
therefore it makes sense to extend our considerations to a
general metrized line bundle $(L,\rho)$. Of course, whenever
referring to a conformal metric or to the variational problem for
the Liouville equation, it will be be assumed that $L=T_X$ and
that $g_{ij}= z'_{ij}$. Upon choosing a cover $\cover{U}_X$,
which for now we assume to be a \cech\ cover by open sets, the
pair $(L,\rho)$ is described in terms of the data expounded in
section~\ref{sec:line-bun}. Our starting point will be the
$0$-cochain
\begin{equation}
  \label{eq:48}
  \omega^{0}_i[\log\rho_i]
  =\onehalf\,\del \log \rho_i \wedge \delb \log \rho_i
  = -\onehalf d \log \rho_i \wedge \onehalf d^c \log \rho_i
\end{equation}
with values in $\she[2]{X}(1)(U_i)$.
\begin{remark}
  \label{rem:1}
  A generalization for eq.~\eqref{eq:48} would be to consider a
  \emph{pair} of metrized line bundles $(L,\rho)$ and $(L',\rho')$, and
  then the analog of~\eqref{eq:48} would be
  \begin{equation}
    \label{eq:49}
    \omega^{0}_i[\log\rho_i,\log\rho'_i] = \onehalf \Bigl(
    -\onehalf d \log \rho_i \wedge \onehalf d^c \log \rho'_i
    +\onehalf d^c \log \rho_i \wedge \onehalf d \log \rho'_i
    \Bigr)
  \end{equation}
  Note, however, that the expressions are quadratic. Moreover,
  $L\otimes L'$ has metric $\rho\rho'$, so that there is a
  natural ``polarization identity''
  \begin{equation}
    \label{eq:50}
    \omega[\log\rho,\log\rho'] =
    \onefourth \omega[\log \rho\rho']
    - \onefourth \omega[\log \frac{\rho}{\rho'}]
  \end{equation}
  where we have omitted the indexes for simplicity of notation.  We
  shall comment later on the significance of eq.~\eqref{eq:50}.
\end{remark}

\subsubsection{Computation}
\label{sec:computation}

Let us extend~\eqref{eq:48} to a degree $2$ cocycle in the total
simple complex associated to the double complex
\begin{math}
  \Check{C}^\bullet(\cover{U}_X,\she{X}(1))
\end{math}
of \cech\ cochains with values in the de~Rham complex of imaginary
smooth forms. This is accomplished in the usual fashion (see
e.g.~\cite{bott1982}) by finding a $1$-cochain of $1$-forms
$\omega^{1}_{ij}[\log \rho]$ on $U_{ij}$ and a $2$-cochain of $0$-forms
$\omega^{2}_{ijk}[\log \rho]$ on $U_{ijk}$ such that the relations
\begin{equation*}
  \deltacheck \omega^{0} = -d \omega^{1}\,, \qquad
  \deltacheck \omega^{1} = d \omega^{2}\,, \qquad
  \deltacheck \omega^{2} = 0
\end{equation*}
are satisfied. Of course, the remaining one, namely $d\omega =0$
is automatically satisfied for dimensional reasons. It turns out
that to a great extent these relations are explicitly computable
without further assumptions, such as that the cover $\cover{U}_X$
be good.  The needed calculations are fairly standard, and they
are presented in great detail in ref.~\cite{math.CV/0204318}, so
we shall be brief. (The observation in \cite{lang:arakelov} that
on a Riemann surface for two smooth functions $f$ and $g$ one has
\begin{math}
  df \wedge d^c g = dg\wedge d^c f
\end{math}
is useful in carrying out the calculations.) The first two steps
are as follows.

First, one has:
\begin{equation}
  \label{eq:51}
  \begin{gathered}
    \omega^{0}_j[\rho] - \omega^{0}_i[\rho]
    = -d\,\omega^{1}_{ij} [\rho] \\
    \omega^{1}_{ij} [\rho]
    = \onehalf\log\rho_i\;d^c \log\,\abs{g_{ij}}
    +d^c \log\,\abs{g_{ij}}\; \onehalf \log \rho_j
  \end{gathered}
\end{equation}
The next step yields:
\begin{equation*}
  \begin{aligned}
    \deltacheck\bigl( \omega^{1}[\log\rho] \bigr)_{ijk} &=
    \omega^{1}_{jk} [\rho] -\omega^{1}_{ik} [\rho]
    +\omega^{1}_{ij} [\rho] \\
    &= \log\,\abs{g_{ij}} \; d^c \log\,\abs{g_{jk}}
    -d^c \log\,\abs{g_{ij}} \; \log\,\abs{g_{jk}}
  \end{aligned}
\end{equation*}
and notice that
\begin{math}
  d^c \log\,\abs{g_{ij}} = \pi_1 d\,\log \,g_{ij}\,,
\end{math}
and
\begin{math}
  \pi_{p+q-1}(a\wedge b) = \pi_p(a)\wedge \pi_{q-1}(b) +
  \pi_{p-1}(a)\wedge \pi_q(b)\,,
\end{math}
so we have
\begin{equation}
  \label{eq:52}
  \deltacheck\bigl( \omega^{1}[\log\rho] \bigr)_{ijk} =
  \pi_1\bigl( \log\,g_{ij} \; d\,\log\, g_{jk} \bigr)
  -d \bigl( \pi_1(\log\,g_{ij})\; \log\,\abs{g_{jk}} \bigr) \,.
\end{equation}
Observe that now the problem of continuing the descent becomes
independent of the chosen metric $\rho$. The most direct way of
proceeding is the following. If we assume the cover $\cover{U}_X$ to be
acyclic for the de~Rham complex $\she{X}(1)$, then there exists
$\omega^{2}_{ijk} \in \she[0]{X}(1)(U_{ijk})$ such that
$d\omega^{2}_{ijk} = \deltacheck (\omega^{1}[\log\rho])_{ijk}$.
Furthermore, consistency on a quadruple intersection requires that
$\deltacheck\omega^{2}$ be a $3$-cocycle with values in $\RR(1)_X$.
Since $H^3(X,\RR(1))=0$, this cocycle must be a coboundary, therefore,
up to readjusting the constants, there exists a choice of
$\omega^{2}_{ijk}$ such that $\deltacheck\omega^{2} = 0$, and
furthermore, $\omega^{2}_{ijk}$ does not depend on the metric structure.

\subsubsection{Solution to the variational problem}
\label{sec:solut-vari-probl}

The previous preliminary calculation is sufficient from the point of
view of finding the extrema. To this effect, we set $L = T_X$, for
$\rho$ a conformal metric in $\cm{X}$. Notice that the space of
conformal metrics on $X$ is affine over $C^\infty(X,\RR) \equiv
\mathcal{E}^0(X)$: if $ds^2$ and ${ds'}^2$ are two conformal metrics
with local expressions $\rho_i\,\abs{dz_i}^2$ and
$\rho'_i\,\abs{dz_i}^2$ respectively, then there exists $\sigma \in
C^\infty(X,\RR)$ such that $\log\rho'_i =\sigma\vert_{U_i} +
\log\rho_i$. The change from $\Omega[\log\rho]$ to $\Omega[\log\rho
+\sigma]$ can be exactly computed thanks to the fact that the last step
in the determination of $\Omega[\log\rho]$ is independent of $\rho$ and
the quadratic character of~\eqref{eq:44}. Indeed we have:
\begin{proposition}
  \label{prop:2}
  \begin{equation*}
    S[\log\rho +\sigma] -S[\log\rho] =
    \frac{1}{\tate} \int_X \bigl(
    \onehalf \del\sigma \wedge \delb\sigma
    +\sigma\, c_1(\rho) -\sqrt{-1}(e^\sigma -1)
    \omega_\rho\bigr)
  \end{equation*}
\end{proposition}
\begin{proof}
  The change in $\omega^{0}_i[\log \rho]$ is computed as
  \begin{equation*}
    \omega^{0}_i[\log \rho + \sigma] -\omega^{0}_i[\log \rho] =
    -\onehalf d\sigma_i \wedge \onehalf d^c\sigma_i
    + \sigma_i\, \onehalf dd^c \log \rho_i
    - \onehalf d(\sigma_i \, d^c \log\rho_i)\,,
  \end{equation*}
  where we set $\sigma_i \equiv \sigma\rvert_{U_i}$. Note that
  the first two terms on the right hand side are globally
  well-defined $2$-forms. On the other hand
  \begin{equation*}
    \omega^{1}_{ij}[\log \rho +\sigma]
    -\omega^{1}_{ij}[\log\rho] =
     \sigma_i\,d^c \log\,\abs{g_{ij}}\,.
  \end{equation*}
  Letting $\chi_i = \sigma_i \, d^c \log\rho_i$, we see that
  \begin{equation*}
    \Omega[\log \rho +\sigma] - \Omega[\log \rho] =
    -\onehalf d\sigma \wedge \onehalf d^c\sigma
    + \sigma\, \onehalf dd^c \log \rho -D\chi\,,
  \end{equation*}
  and taking the area terms into account, establishes the
  formula.
\end{proof}
As a consequence, we obtain
\begin{corollary}
  \label{cor:1}
  The Liouville equation~\eqref{eq:42} is the Euler-Lagrange
  equation for the Liouville functional~\eqref{eq:47} introduced
  in definition~\ref{def:2}. The critical point is
  non-degenerate.
\end{corollary}
\begin{proof}
  Replacing $\sigma$ with $t\,\sigma$, $t\in\RR$, in the previous
  proposition we find the infinitesimal change in $\Omega[\log
  \rho]$ to be
  \begin{equation*}
    \left.\frac{d}{dt}\right\rvert_{t=0}
    S[\log \rho + t\sigma] = -\frac{1}{\tate}\, \int_X \sigma
    \bigl( c_1(\rho) - \sqrt{-1}\,\omega_\rho \bigr)\,,
  \end{equation*}
  and it follows that $S[\log\rho]$ has an extremum if and only
  if the Liouville equation is satisfied. Non degeneracy follows
  from the quadratic part in the exact change formula in
  proposition~\ref{prop:2}.
\end{proof}
\begin{remark}
  The fact that the change in the cocycle is given by a pure $2$-form
  term up to total coboundary can also be analyzed in terms of gluing
  properties of variational bicomplexes, cf.  Theorem 1.2 in ref.\
  \cite{MR1908413} and the proof of Theorem 1 in ref.\ \cite{aldtak2000}.
  From this perspective, Corollary~\ref{cor:1} is a direct consequence
  of the affine structure of the space $\cm{X}$ of conformal metrics.
\end{remark}

\subsection{A cup product}
\label{sec:cup-product}
Formula~\eqref{eq:52} can be handled in a more geometric fashion
as follows. From sec.~\ref{sec:remarks-dilogarithm} we can
rewrite~\eqref{eq:52} as
\begin{equation}
  \label{eq:53}
  \deltacheck\bigl( \omega^{1}[\log\rho] \bigr)_{ijk} =
   -r_2(g_{ij}, g_{jk})\,,
\end{equation}
and we have the collection of tame symbols
$\tame{g_{ij}}{g_{jk}}$ associated any triple intersection
$U_{ijk}$ in the cover $\cover{U}_X$, \cite{brymcl:deg4_II}.
These symbols glue to form a global symbol $\tame{L}{L}$. As a
cohomology class on a curve $X$, however, $\tame{L}{L}$ will be
zero (that is, there is a global object in the associated
$2$-stack, cf. \cite{bry:quillen}), so that it will be possible
to choose local functions $L_{ijk}$ such that
\begin{equation*}
  d L_{ijk} = -\log\,g_{ij} \; d\,\log\, g_{jk}\,,
\end{equation*}
as explained in detail in \cite{aldtak2000}. Note that we still need the
cover $\cover{U}_X$ to be fine enough. Moreover, cf. loc.~cit., the
collection $L_{ijk}$ can be chosen in a way that
\begin{equation*}
  \deltacheck L _{ijkl} = -c_{ijk}  \log \, g_{kl} + n_{ijkl}\,,
\end{equation*}
where $n_{ijkl} \in \ZZ(2)$. Therefore
from~\eqref{eq:21},~\eqref{eq:52} and~\eqref{eq:53}, we can set
\begin{equation}
  \label{eq:54}
  \omega^{2}_{ijk} = 
  -\bigl( \pi_1(\log\,g_{ij})\; \log\,\abs{g_{jk}} \bigr)
  -\pi_1 L_{ijk}\,,
\end{equation}
namely according to sec.~\ref{sec:remarks-dilogarithm} we have set
$\omega^{2}_{ijk} = -\bwli (g_{ij},g_{jk})$. Now, the last compatibility
condition is satisfied, indeed we have:
\begin{align*}
  \deltacheck \omega^{2}_{ijkl} &= 
  -c_{ijk}\, \log \,\abs{g_{kl}} -\pi_1 (\deltacheck L_{ijkl}) \\
  &= -c_{ijk}\, \log \,\abs{g_{kl}}
  +\pi_1 (c_{ijk} \log g_{kl} ) \\
  &= 0\,.
\end{align*}
As a result, we obtain a $2$-cocycle in the the total complex associated
to $\Check{C}^\bullet (\cover{U}_X,\she{X}(1))$ as before, with a more
geometric interpretation of $\omega^{2}_{ijk}$ in terms of the
trivialization of the symbol $\tame{L}{L}$. Moreover, notice that if $X$
is a curve then $\she{X}(1)[1] \qi \deltilde{\RR}{2}$ thus we may
interpret the class so determined by $\Omega[\log \rho]=
\omega^{0}[\log\rho] +\omega^{1}[\log\rho] +\omega^{2}$ as a (degree
$3$) class in $\delH[3]{X}{\RR}{2}$.

\subsection{Two line bundles}
\label{sec:two-line-bundles}

For a pair $(L,\rho)$, $(L',\rho')$ of metrized line bundles we can
complete~\eqref{eq:49} to a cocycle $\Omega[\log\rho,\log\rho']$ via an
analogous procedure to the one presented in
sects.~\ref{sec:Direct-constr-Liouv} and~\ref{sec:cup-product}.  The
relevant calculations being entirely similar, we limit ourselves to
quoting the relevant expressions. Starting from~\eqref{eq:49}, which we
rewrite in the form
\begin{equation}
  \label{eq:55}
  \omega^{0}_i[\log\rho_i,\log\rho'_i] =
  \onehalf d^c \log \rho_i \wedge \onehalf d \log \rho'_i\,,
\end{equation}
the corresponding expression for the degree $(1,1)$ term is:
\begin{equation}
  \label{eq:56}
  \omega^{1}_{ij}[\log\rho,\log\rho'] =
  \onehalf \log\rho_i\;d^c \log\,\abs{g'_{ij}} 
  +\onehalf d^c \log\,\abs{g_{ij}}\; \log \rho'_j\,.
\end{equation}
Computing the \cech\ coboundary we find:
\begin{equation}
  \label{eq:57}
  \begin{aligned}
    \bigl(\deltacheck \omega^{1}[\log\rho,\log\rho']\bigr)_{ijk} &=
    \log\,\abs{g_{ij}} \; d^c \log\,\abs{g'_{jk}}
    -d^c \log\,\abs{g_{ij}} \; \log\,\abs{g'_{jk}} \\
    &= -r_2(g_{ij},g'_{jk})\,,
  \end{aligned}
\end{equation}
from which $\omega^{2}_{ijk}$ (now independent of $\rho$ and $\rho'$)
can be obtained as $-\bwli (g_{ij},g'_{jk})$ by looking at a collection
$L_{ijk}$ such that
\begin{equation*}
  d L_{ijk} = -\log\,g_{ij} \; d\,\log\, g'_{jk}\,,
\end{equation*}
from the triviality of the symbol $\tame{L}{L'}$.

\subsection{Additional remarks on the Liouville functional and the
  Bloch-Wigner dilogarithm}
\label{sec:liouv-funct-bloch}

We wish to compare the previous constructions to those of ref.\ 
\cite{math.CV/0204318}. To this end we need to specifically
consider the case of a cover of $X$ provided by a Kleinian
uniformization. (The reader should consult loc.\ cit.\ for
reference and complete details.) The comparison offers a better
perspective on the absence of cohomological obstructions in the
calculations of ref.\ \cite{math.CV/0204318}, and on the
relations with three-dimensional hyperbolic geometry.

Let $\Gamma$ be a purely loxodromic second kind Kleinian group
satisfying all the conventions spelled out in
Appendix~\ref{sec:kleinian-groups}, to which we refer for the
notation.  Let $U_\Gamma\subset \pione$ be the region of
discontinuity, and $X=U_\Gamma/\Gamma$ the resulting (possibly
disconnected) Riemann surface.

A conformal metric $\rho$ on $X=X_1\sqcup \dots \sqcup X_n$
appears as an automorphic function on $U_\Gamma$:
\begin{equation}
  \label{eq:62}
  \onehalf \log \rho - \onehalf \log \rho\circ \gamma
  = \log\, \abs{\gamma'}\,,\quad \gamma\in \Gamma\,.
\end{equation}
Eq.~\eqref{eq:62} is the direct translation of \eqref{eq:29}
following the principles of \cite{aldtak2000}.  Accordingly, the
first two terms of the Liouville cocycle computed by applying the
procedure explained in sect.~\ref{sec:Conf-metr-Liouv} are:
\begin{gather}
  \omega^{0}[\log\rho]
  = -\onehalf d \log \rho \wedge \onehalf d^c \log \rho \label{eq:63} \\
  \omega^{1}_{\gamma} [\log\rho]
  = \bigl(\onehalf\log\rho +\onehalf \log \rho
  \circ \gamma \bigr)\; d^c \log\,\abs{\gamma'}\label{eq:64}
\end{gather}
and computing the coboundary of~\eqref{eq:64} according to the
prescription in Appendix~\ref{sec:kleinian-groups} yields
\begin{equation}
  \label{eq:65}
  \deltacheck\bigl( \omega^{1}[\log\rho]
  \bigr)_{\gamma_1,\gamma_2}
  =-r_2((\gamma_1\gamma_2)', \gamma'_2)\,,
\end{equation}
where $r_2$ has been introduced in eqs.~\eqref{eq:18}
and~\eqref{eq:19}. Note that, as in sect.~\ref{sec:cup-product},
the coboundary of $\omega^{1}$ is a cup product in real Deligne
cohomology:
\begin{equation*}
  r_2((\gamma_1\gamma_2)', \gamma'_2)
  = (d\log (\gamma_1\gamma_2)', \log\, \abs{(\gamma_1\gamma_2)'})
  \cup (d\log \gamma'_2, \log\, \abs{\gamma'_2})\,.
\end{equation*}
where the two classes are associated to the rational functions
$(\gamma_1\gamma_2)'$ and $\gamma'_2$, respectively (cf.\ 
sect.~\ref{sec:remarks-dilogarithm}). Hence, we can work with the
double complex of group cochains on $\Gamma$ with values in the
real Deligne complex $\deltilde{\RR}{2}\qi \she{U_\Gamma}(1)$ on
the region of discontinuity $U_\Gamma$. Since from this point on
only rational functions with singularities at certain prescribed
points will appear, following ref.\ \cite{math.AG/0003086} we
will consider the Deligne complex on the generic point
$\eta_{\pione}$ of $\pione$.

For any two elements $\gamma_i, \gamma_j \in \Gamma$ define
$T_{ij}\in \PSL_2(\CC)$ by
\begin{equation}
  \label{eq:66}
  z \longmapsto T_{ij} (z) = \CR{z}{z_{ij}}{z_j}{\infty}
  = \frac{z-z_{ij}}{z-z_j}\,.
\end{equation}
where $z_j = \gamma_j^{-1}(\infty)$ and $z_{ij} =
\gamma_j^{-1}(z_i)$. Following ref.\ \cite{math.CV/0204318}, we
introduce the $1$-cochain on $\Gamma$ with values in
$\she[1]{\pione}(1)(\eta_{\pione})$:
\begin{equation}
  \label{eq:68}
  \kappa_\gamma = \log\,\abs{c_\gamma} \, \pi_1 d\log \gamma'\,.
\end{equation}
We have:
\begin{lemma}
  \label{lemma:kappa}
  \begin{math}
    r_2((\gamma_1\gamma_2)', \gamma'_2) =
    4\, d \bigl( \bwli \circ T_{12} \bigr)
    +\deltacheck\kappa_{\gamma_1,\gamma_2}\,,
  \end{math}
  where $\bwli = \sqrt{-1} \bwd$, and $\bwd$ is the standard
  Bloch-Wigner dilogarithm, cf.\ \ref{sec:remarks-dilogarithm}.
\end{lemma}
\begin{proof}
  A straightforward calculation exploiting relation~\eqref{eq:60}.
\end{proof}
Since obviously $d\kappa_\gamma=0$, we can redefine
\begin{math}
  \omega^{1}_\gamma \to \omega^{1}_\gamma + \kappa_\gamma\,,
\end{math}
so that using the lemma, from eq.~\eqref{eq:65} we have:
\begin{equation}
  \label{eq:69}
  \deltacheck\bigl( \omega^{1}[\log\rho]
  \bigr)_{\gamma_1,\gamma_2} 
  = -4\, d\,\bigl(\bwli \circ T_{12} \bigr)\,.
\end{equation}
For convenience of notation, let us temporarily set
\begin{math}
  \Hat\omega^{2}_{\gamma_1,\gamma_2} := -4\, \bigl(\bwli \circ
  T_{12} \bigr)\,.
\end{math}
We then have the following
\begin{lemma}
  \label{lem:3}
  \begin{equation}
    \label{eq:71}
    D(\omega^{0}+\omega^{1}+\Hat\omega^{2}) =
    -4 \,\bwli (
    \CR{\infty}%
    {\gamma_1(\infty)}{\gamma_1\gamma_2(\infty)}%
    {\gamma_1\gamma_2\gamma_3(\infty)}) \,. 
  \end{equation}
\end{lemma}
\begin{proof}
  By construction,
  \begin{math}
    D(\omega^{0}+\omega^{1}+\Hat\omega^{2}) =
    \deltacheck\Hat\omega^{2}\,,
  \end{math}
  and the latter \cech\ coboundary is computed as
  \begin{equation*}
    \begin{split}
      \deltacheck\Hat\omega^{2}_{\gamma_1,\gamma_2, \gamma_3}
      &= -4 \Bigl(\bwli (\CR{z}{z_{23}}{z_3}{\infty})
      -\bwli (\CR{z}{z_{123}}{z_3}{\infty})\\
      &\qquad +\bwli (\CR{z}{z_{123}}{z_{23}}{\infty})
      -\bwli (\CR{\gamma_3(z)}{z_{12}}{z_2}{\infty}) \Bigr)
    \end{split}
  \end{equation*}
  In the last term in the previous relation we have
\begin{math}
  \CR{\gamma_3(z)}{z_{12}}{z_2}{\infty}
  =\CR{z}{z_{123}}{z_{23}}{z_3}\,.
\end{math}
The Bloch-Wigner dilogarithm satisfies the $5$-term
relation~\cite{MR2002g:52013,MR94k:19002}:
\begin{equation}
  \label{eq:70}
  \sum_{i=0}^4 (-1)^i \bwd ([a_0:\dots :\Hat{a_i}:\dots :a_4]) = 0\,,
\end{equation}
where $a_0,\dots, a_4 \in \pione$ and the hat sign denotes
omission. As a consequence we have:
\begin{equation*}
  \deltacheck\Hat\omega^{2}_{\gamma_1,\gamma_2, \gamma_3} =
  -4 \,\bwli (\CR{z_{123}}{z_{23}}{z_3}{\infty}) \,,
\end{equation*}
and again by the invariance of the cross ratio, we
obtain~\eqref{eq:71}.
\end{proof}
\begin{remark}
  \begin{math}
    \bwd (\CR{a}{b}{c}{d} )
  \end{math}
  is the hyperbolic volume of the ideal hyperbolic tetrahedron
  with vertices at the points
  \begin{math}
    a,b,c,d \in \pione\,,
  \end{math}
  see, e.g.\ refs.\ \cite{MR84b:53062b,MR88k:57032,MR94k:19002}.
\end{remark}
It follows from the five-term relation~\eqref{eq:70} that the
right hand side of eq.~\eqref{eq:71} defines an $\RR(1)$-valued
$3$-cocycle on $\Gamma$. Moreover, this cocycle is already defined
on $\PSL_2(\CC)$, where its class is known to generate the
Eilenber-Mac~Lane cohomology group
\begin{math}
  H^3(\PSL_2(\CC),\RR(1))\,,
\end{math}
\cite{MR2001i:11082,MR88k:57032,MR94k:19002}.  It is also known
that up to a factor $2$ it agrees with the imaginary part of the
second Cheeger-Simons universal secondary class $\Hat C_2$.

Thus the complete Liouville cocycle
\begin{math}
  \Omega = \omega^{0}+\omega^{1}+\omega^{2}
\end{math}
subordinated to the cover $U_\Gamma \to X$ is found as follows.
The pullback of $\Hat C_2$ along the inclusion map
\begin{math}
  \Gamma \hookrightarrow \PSL_2(\CC)
\end{math}
is zero, since ($\Gamma$ being of second kind) the
$3$-manifold $M_\Gamma = \HH^3/\Gamma$ is non-compact, and
$H^\bullet(\Gamma,\RR(1)) \cong H^\bullet(M_\Gamma,\RR(1))$. It
follows that the restriction of the cocycle given by the
Bloch-Wigner dilogarithm to $\Gamma$ must be a coboundary, hence
there exists a group $2$-cochain $c$ on $\Gamma$ with values in
$\RR(1)$ such that
\begin{equation}\label{eq:67}
 4 \,\bwli (
  \CR{\infty}%
  {\gamma_1(\infty)}{\gamma_1\gamma_2(\infty)}%
  {\gamma_1\gamma_2\gamma_3(\infty)})
  = (\deltacheck c)_{\gamma_1,\gamma_2,\gamma_3}\,.
\end{equation}
It follows that the cochain $c$ provides the necessary
``integration constants,'' namely the required $2$-cochain on
$\Gamma$ with values in $\she[0]{U_\Gamma}$ to complete the
Liouville cocycle is
\begin{math}
  \omega^{2}_{\gamma_1,\gamma_2}
  = -4\, \bigl(\bwli \circ
  T_{12} \bigr)
  + c_{\gamma_1,\gamma_2}\,.
\end{math}

\section{Conformal metrics and hermitian holomorphic cohomology}
\label{sec:conf-metr-herm}

In sect.~\ref{sec:Conf-metr-Liouv} we presented a construction of a
degree $2$, $\RR(1)$-valued class corresponding to a conformal metric
$\rho\in \cm{X}$, represented by the cocycle $\Omega[\log\rho]$.
Supplemented by the area of $X$ computed with respect to $\rho$, it
provides a global functional for the variational problem associated to
the Liouville equation~\eqref{eq:42}. We now show that it coincides
with the square of the class of $(T_X,\rho)$ in hermitian holomorphic
Deligne cohomology introduced in sect.~\ref{sec:herm-holom-deligne}.
Moreover, we show that this equality holds at the cocycle level. More
generally, without considering the area term, we show that the cup
product of $(L,\rho)$ and $(L',\rho')$ in hermitian holomorphic Deligne
cohomology coincides with the class of $\Omega[\log\rho,\log\rho']$, and
again the equality in fact holds at the cocycle level.

\subsection{Comparison on a curve}
\label{sec:comparison-curve}

Let $X$ be a compact Riemann surface. From the results of
sect.~\ref{sec:Cup-prod-herm}, the cup product of the classes of
$(L,\rho)$ and $(L',\rho')$ in hermitian holomorphic Deligne cohomology
yields a class in $\dhhH[4]{X}{2}$, and on a curve we only capture the
$2$-dimensional part of this class. Indeed, in the exact
sequence~\eqref{eq:25}, the cohomology class corresponding to the symbol
$\tame{L}{L'}\in \delH[4]{X}{\ZZ}{2}$ is zero, and
$A^{(2,2)}(X)_{\RR(2)}$ is also zero for obvious dimensional reasons, so
we have:
\begin{equation*}
  \dotsm \longrightarrow \delH[3]{X}{\RR}{2}
  \longrightarrow
  \dhhH[4]{X}{2} \longrightarrow 0\,.
\end{equation*}
It follows that the class $[L,\rho]\cup [L',\rho']$ must come from an
element in $\delH[3]{X}{\RR}{2}$. As already remarked, on a curve we
have $\deltilde{\RR}{2} \qi \she{X}(1)[-1]$, thus
$\delH[3]{X}{\RR}{2}\cong \HHH^3(X,\she{X}(1)[-1]) \cong
\HHH^2(X,\she{X}(1)) \cong H^2(X,\RR(1))$, in agreement with the
calculations performed in sect.~\ref{sec:Conf-metr-Liouv}.

More in detail, in complex dimension $1$ the second hermitian
holomorphic Deligne complex $\dhh{1}$ simplifies considerably and
diagram~\eqref{eq:39} reduces to
\begin{equation}
  \label{eq:72}
  \begin{CD}
    \ZZ(2) @>-\imath>> \sho{X} @>-d>> \shomega[1]{X} \\
    & & @VV\pi_1V @VV\pi_1V \\
    & & \she[0]{X}(1) @>-d>> \she[1]{X}(1) @>-d>> \she[2]{X}(1)
  \end{CD}
\end{equation}
so that $\dhh{2}$ becomes just the cone of the morphism
$\deligne{\ZZ}{2} \overset{\pi_1}{\longrightarrow}
\she{X}(1)[-1]$. In other words, on a curve $X$ we have that
$\dhh{2}$ is given by the complex
\begin{equation}
  \label{eq:73}
  \begin{CD}
    \ZZ(2)_X @>{-\imath}>>
    \sho{X}  @>{(-d,-\pi_1)}>>
    \shomega[1]{X}\oplus \she[0]{X}(1) @>{-\pi_1+d}>>
    \she[1]{X}(1) @>{d}>>
    \she[2]{X}(1)\,,
  \end{CD}
\end{equation}
where the differentials have been written explicitly. We can see
the complex $\she{X}$ appears as a subcomplex in~\eqref{eq:73}
and the shift of two positions to the right clearly accounts for
the cohomology degree shift from $2$ to $4$.

Our main result is the following comparison
\begin{theorem}
  \label{thm:2}
  Let $X$ be a compact Riemann surface of genus $g>1$. Let $(L,\rho)$
  and $(L',\rho')$ be two hermitian holomorphic line bundles. The class
  of $[L,\rho]\cup [L',\rho']$ in $\dhhH[4]{X}{2} \cong H^2(X,\RR(1))$
  coincides with the one represented by the cocycle $\Omega[\log
  \rho,\log \rho']$ constructed in section~\ref{sec:Conf-metr-Liouv}.
\end{theorem}
\begin{proof}
  We have observed above that $\dhhH[4]{X}{2} \cong
  H^2(X,\RR(1))$, and by construction the class of
  \begin{math}
    \Omega[\log\rho,\log\rho']
  \end{math}
  is in $\HHH^2(X,\she{X}(1)) \cong H^2(X,\RR(1))$. Note that for $X$
  connected they must coincide up to a proportionality factor, since
  $H^2(X,\RR(1))\cong \RR(1)$ in this case. In general, we compute the
  proportionality factor from the explicit cocycles from
  sects.~\ref{sec:Cup-prod-herm} and~\ref{sec:Direct-constr-Liouv}
  to~\ref{sec:two-line-bundles}. (Since $\Omega[\log\rho,\log\rho']$ is
  computed under suitable acyclicity assumptions on the cover, so we
  will use such a cover to establish the comparison.)
  
  Let us assume $L$ and $L'$ and their respective hermitian metric
  structures are represented by cocycles of type~\eqref{eq:38} with
  respect to the chosen cover $\cover{U}_X$. Specializing the general
  expression in table~\ref{tab:big-cup} in sect.~\ref{sec:Cup-prod-herm}
  to the case at hand we obtain, with reference to~\eqref{eq:73}:
  \begin{equation}
    \label{eq:74}
    \renewcommand{\arraystretch}{1.3}
    \begin{array}{cc}
      (0,4) & c_{ijk}c'_{klm} \\ \hline
      (1,3) & -c_{ijk} \log \, g'_{kl} \\ \hline
      (2,2) & -\log \, g_{ij} \, d\log\,g'_{jk} \oplus
      -\pi_1\log\,g_{ij}\, \log\,\abs{g'_{jk}} \\ \hline
      (3,1) &
        (1-\alpha) \bigl( d^c\log\,\abs{g_{ij}}\; \onehalf \log\rho'_j
        -\log\,\abs{g_{ij}}\; \onehalf d^c \log\,\rho'_j
        \bigr) 
        +\alpha \bigl( -\onehalf d^c \log\rho_i\;\log\,\abs{g'_{ij}}
        +\onehalf \log\rho_i\; d^c\log\,\abs{g'_{ij}} \bigr) \\ \hline
      (4,0) &
      \alpha\, \delb\del \log\rho_i\;
        \onehalf\log\rho'_i +(1-\alpha)\, \onehalf \log\rho_i \;
        \delb\del\log\rho'_i
    \end{array}
  \end{equation}
  where we have followed the convention explained in the introduction
  for the bidegrees. Let us denote by $\theta^{i}$ the term of bidegree
  $(4-i,i)$ in~\eqref{eq:74} and by $\Theta$ the total cocycle. (For
  simplicity, we suppress $\rho$ and $\rho'$ from the notation.) A
  direct calculation shows that
  \begin{equation*}
    \theta^{0}_i = \omega^{0}_i + d\lambda^{0}_i\,,\qquad
    \theta^{1}_{ij} = \omega^{1}_{ij} -\deltacheck\lambda^{0}_{ij}\,,
  \end{equation*}
  where $\omega^{0}_i$ is given by eq.~\eqref{eq:55}, $\omega^{1}$
  is given by eq.~\eqref{eq:56}, and
  \begin{equation*}
    \lambda^0_i = \alpha\,
    \onehalf d^c \log\rho_i \, \log\rho'_i
    +(1-\alpha)\, \onehalf \log\rho_i \, d^c \log\rho'_i\,.
  \end{equation*}
  Note that at this point we could simply \emph{define} $\Omega = \Theta
  -D\lambda^0$. Furthermore, note that $\Omega$ does not explicitly
  depend on the parameter $\alpha$ from the \bei\ product~\eqref{eq:5}.
  To finish the comparison, let us assume the cover $\cover{U}_X$ allows
  us to find a collection $L_{ijk}\in \sho{X}(U_{ijk})$ such that
  \begin{equation}
    \label{eq:75}
    \begin{aligned}
      d L_{ijk} &= -\log\,g_{ij} \; d\,\log\, g'_{jk}\\
      \deltacheck L _{ijkl} &= -c_{ijk}  \log \, g'_{kl} + n_{ijkl}\,,
    \end{aligned}
  \end{equation}
  as in sect.~\ref{sec:cup-product}. In this way we have 
  \begin{math}
    \Omega =
    (\omega^{0},\omega^{1}, \omega^{2})
  \end{math}
  with $\omega^{2}_{ijk} = -\bwli (g_{ij},g'_{jk})$. $\Omega$ is a
  cocycle of total degree $2$ in $\Tot \Check{C}^\bullet
  (\cover{U}_X,\she{X})$, and it injects (via the exact sequence of the
  cone) into $\Tot\Check{C}^\bullet (\cover{U}_X,\dhh{2})$ as
  \begin{equation*}
    (\omega^{0},\omega^{1}, 0\oplus\omega^{2})\,.
  \end{equation*}
  Then via equations~\eqref{eq:75} it is easily seen that
  \begin{equation*}
    \theta^{2}_{ijk}=d_{\dhh{2}} (-L_{ijk})
    + 0\oplus \omega^{2}_{ijk}\,,
  \end{equation*}
  where $d_{\dhh{2}}$ is the differential in $\dhh{2}$, and therefore
  \begin{equation*}
    \Theta = (\omega^{0},\omega^{1}, 0\oplus \omega^{2})
    +D\lambda^{0} +D(-L,n)\,,
  \end{equation*}
  where we have put $D=d_{\dhh{2}}\pm \deltacheck$ for the total
  differential. Thus the two cocycles constructed via the direct method
  of sect.~\ref{sec:Direct-constr-Liouv} and the cup product of metrized
  bundles define the same class. By direct comparison, the
  proportionality factor is $1$.
\end{proof}
In light of the previous theorem, the polarization identity in
Remark~\ref{rem:1} is now easily explained. Using
\begin{math}
  [L\otimes L',\rho\rho'] = [L,\rho] + [L',\rho']
\end{math}
and
\begin{math}
  [L\otimes {L'}^\vee,\rho/\rho'] = [L,\rho] - [L',\rho'] \,,
\end{math}
and the (graded) commutativity of the cup product
\begin{equation*}
    \dhhH[2]{X}{1} \otimes \dhhH[2]{X}{1}
  \overset{\cup}{\longrightarrow}
  \dhhH[4]{X}{2}\,,
\end{equation*}
we obtain the polarization identity
\begin{equation*}
  4 [L,\rho]\cup [L',\rho'] = [L\otimes L',\rho\rho']^2
  -[L\otimes {L'}^\vee,\rho/\rho']^2\,,
\end{equation*}
where the squares in the right hand side refer to $\cup$. A polarization
identity at the level of representative cocycles, and hence the one in
Remark~\ref{rem:1}, follow by applying Thm.~\ref{thm:2} to the latter
identity.

By choosing $L=L'=T_X$, the holomorphic tangent
bundle of $X$, we immediately obtain:
\begin{corollary}
  \label{corollary:1}
  Let $\rho\in\cm{X}$ be a conformal metric. The Liouville functional
  without area term~\eqref{eq:46} is given by the (evaluation of) the
  square $[T_X,\rho]\cup [T_X,\rho]$ with respect to the cup product in
  hermitian holomorphic Deligne cohomology. The full-fledged Liouville
  functional is obtained by adding the area term
  \begin{math}
    \frac{1}{2\pi}\int_X \omega_\rho
  \end{math}
  to~\eqref{eq:46}.
\end{corollary}
\begin{remark}
  Due to the specific form of the differential in the complex $\dhh{2}$
  the descent equations are explicit and close automatically. Therefore
  the cocycle resulting from the calculation of the cup product
  sidesteps the problem of the explicit calculation of the last term,
  unlike the more direct version from sect.~\ref{sec:Conf-metr-Liouv}.
  Thus, thanks to the explicit character of the calculation, specific
  assumptions on the nature of the cover $\cover{U}_X$ are not required.
\end{remark}
It follows from Thm~\ref{thm:2}, corollary~\ref{corollary:1} and
the previous remark that definition~\ref{def:2} applies to any
(\'etale) cover $\cover{U}_X \to X$. Indeed, proposition~\ref{prop:2}
from sect.~\ref{sec:Direct-constr-Liouv} can be reformulated at the
cocycle level as follows:
\begin{proposition}
  \label{prop:3}
  Let $X$ be a compact, genus $g>1$ Riemann surface and let
  $\cover{U}_X\to X$ be a cover. For a conformal metric $\rho \in
  \cm{X}$ and $\sigma \in C^\infty(X,\RR)$, there is a cocycle
  $\Hat{\Omega}_{\cover{U}_X} [\log\rho]$ solving the variational
  problem for the Liouville equation.
\end{proposition}
\begin{proof}
  If $\rho$ is a conformal metric, let the pair $(T_X,\rho)$ be
  represented, as an hermitian line bundle, by a cocycle $c(T_X,\rho)$
  with respect to the cover $\cover{U}_X$. We set
  \begin{equation*}
    \Omega[\log\rho] = c(T_X,\rho)\cup c(T_X,\rho)\,,
  \end{equation*}
  and a simple calculation starting from eq.~\eqref{eq:74} yields
  \begin{equation*}
    \Omega[\log\rho +\sigma] - \Omega[\log\rho] =
    \sigma\,c_1(\rho) + \onehalf\,\sigma\, \del\delb\sigma +
    D\chi\,,
  \end{equation*}
  where
  \begin{math}
    \chi_i = \onehalf \sigma \, \onehalf d^c\log \rho_i
    -\onehalf d^c\sigma\, \onehalf \log\rho_i\,.
  \end{math}
  Now define
  \begin{equation*}
    \Hat{\Omega}_{\cover{U}_X}[\log\rho] = \Omega[\log\rho]
    -\sqrt{-1}\omega_\rho\,. 
  \end{equation*}
  We see that it yields the formula in Proposition~\ref{prop:2}. In
  particular we have that 
  \begin{equation*}
    \left.\frac{d}{dt}\right\rvert_{t=0}
    \Hat{\Omega}_{\cover{U}_X}[\log\rho +t\sigma]
    \equiv \sigma ( c_1(\rho) -\sqrt{-1}\omega_\rho )\,,
  \end{equation*}
  where $\equiv$ means ``up to total coboundary.''
\end{proof}

\subsection{Determinant of cohomology}
\label{sec:determ-cohom}

Let again $L$ and $L'$ be two holomorphic line bundles with
hermitian metrics $\rho$ and $\rho'$, respectively, on the
compact Riemann surface $X$. Brylinski proves in
\cite{bry:quillen} that the cup product of $L$ and $L'$ in
hermitian holomorphic Deligne cohomology yields the (logarithm)
of the metric $\norm{\cdot}$ on the Deligne pairing $\dual{L}{L'}$
defined in \cite{MR89b:32038}.

It follows, via Thm~\ref{thm:2} and the
isomorphism~\ref{lemma:qi-dhh-bry} between our version of hermitian
holomorphic cohomology and Brylinski's, that the class of
$\Omega[\log\rho,\log\rho']$ is also equal to $\log
\norm{\dual{L}{L'}}$. It is worthwhile to provide a direct proof of this
fact starting from the explicit cocycle given in \eqref{eq:74}.

First, we need to recall a few definitions from \cite{MR89b:32038}. A
complex line $\dual{L}{L'}$ is assigned to the pair $(L, L')$ as
follows.  Let $D$ and $D'$ be divisors on $X$ corresponding to $L$ and
$L'$, and assume they have disjoint supports. Consider two rational
sections, $s$ and $s'$ such that $(s)=D$ and $(s')=D'$. To this datum
one assigns a copy of the complex line generated by the symbol
$\dual{s}{s'}$ subject to the relations:
\begin{equation}
  \label{eq:76}
  \begin{aligned}
    \dual{fs}{s'} &= f(D') \dual{s}{s'}\\
    \dual{s}{gs'} &= g(D) \dual{s}{s'}
  \end{aligned}
\end{equation}
whenever $f$ is a rational function with divisor $(f)$ disjoint from
$D'$, and similarly for $g$. The Weil reciprocity relation
$f(\mathit{div}(g)) = g(\mathit{div}(f))$ (cf. ref.\ \cite{gh:alg_geom})
for two rational functions $f$ and $g$ with disjoint divisors implies
that the relations~\eqref{eq:76} are consistent and the complex line
depends only on the pair $L$, $L'$. When the line bundles are equipped
with hermitian metrics, generically denoted by $\norm{\cdot}$, the
assignment\footnote{We write the square explicitly, whereas the symbol
  $\norm{\cdot}$ used in ref.~\cite{MR89b:32038} denotes the
  \emph{square} of the norm.}
\begin{equation}
  \label{eq:77}
  \log\, \norm{\dual{s}{s'}}^2 = \frac{1}{\tate} 
  \int_X \del\delb\log \norm{s}^2 \,\log\norm{s'}^2
  + \log\norm{s}^2[D'] +\log\norm{s'}^2[D]
\end{equation}
is compatible with the relations~\eqref{eq:76} and defines an hermitian
metric on the complex line $\dual{L}{L'}$. In formula~\eqref{eq:77} the
operator $\del\delb$ is to be computed in the sense of distributions.

Having covered the main definitions, we can now state
\begin{theorem}
  \label{thm:3}
  The cup product of $(L,\rho)$ with $(L',\rho')$ in hermitian
  holomorphic Deligne cohomology corresponds to the logarithm of the
  norm \eqref{eq:77} on the Deligne pairing $\dual{L}{L'}$. The
  proportionality factor is $-\pi\,\sqrt{-1}$.
\end{theorem}
\begin{proof}
Let $D$ and $D'$ be divisors with disjoint support on $X$ corresponding
to $L$ and $L'$, respectively. Using the same technique as in refs.
\cite{bry:quillen,MR89b:32038}, consider two $C^\infty$ positive real
functions $f_1$ and $f_2$ such that $f_1+f_2=1$ and $f_1$ (resp. $f_2$)
vanishes in a neighborhood of the support of $D'$ (resp. $D$). Also, set
$U_1 = X\setminus \mathrm{supp}(D')$ and $U_2 = X\setminus
\mathrm{supp}(D)$. Thus $\{f_1,f_2\}$ is just a partition of unity
subordinated to the cover $\cover{U}_X = \{ U_1, U_2\}$.

The only two terms different from zero in the cocycle $\Theta$
in~\eqref{eq:74} representing the class $[L,\rho]\cup [L',\rho']$ with
respect to this cover are $\theta^{0}_i$ and $\theta^{1}_i$, with the
\cech\ index $i \in \{ 1,2 \}$. Thus the class we are after is
equivalently given by the integral
\begin{equation}
  \label{eq:78}
  \int_X f_1\theta^{0}_1 +f_2\theta^{1}_2
  + df_2 \wedge \theta^{1}_{21} \,,
\end{equation}
which is arrived at by applying in the standard homotopy operator based
on the partition of unity $\{f_1,f_2\}$: from
$\deltacheck\theta^{1}_{12}=0$ we obtain that $\theta^{1}_{12}$ is the
coboundary of the cochain $j\to \sum_{i=1,2} f_i\theta^{1}_{ij}$ and
then we use $\theta^{0}_2-\theta^{0}_2=-d\theta^{1}_{12}$. Observe that
the $2$-form in~\eqref{eq:78} is globally well defined over $X$:
$\theta^{1}_{21}$ is defined only on $U_1\cap U_2$, but $df_2$ has
support on $U_1\cap U_2$, so their wedge product is defined everywhere;
similarly, $f_i\,\theta^{0}_i$ is everywhere defined thanks to the fact
that $f_i$ has support in $U_i$, $i=1,2$.

Consider rational sections $s$ and $s'$ of $L$ and $L'$ such that
$\mathit{div}(s)=D$ and $\mathit{div}(s')=D'$ as above. With respect to
the two-element cover $\cover{U}_X=\{U_1,U_2\}$, the section $s$
corresponds to the pair $\{s_1,s_2\}$, and similarly for
$s'=\{s'_1,s'_2\}$. Since $\mathrm{supp}(D)$ is contained in $U_1$ but
not in $U_2$, and the other way around for $D'$, it follows that $s_2$
and $s'_1$ are actually invertible functions over their respective
domains. Following \cite{lang:arakelov}, we can assume that $s$ and $s'$
are in fact the rational section $\mathbf{1}$, so that $s_2=1$ and
$s'_1=1$, and therefore:
\begin{gather*}
  \abs{s}^2_\rho =
  \begin{cases}
    \log \rho_1 + \log \abs{s_1}^2 & \text{on } U_1\,,\\
    \log \rho_2 & \text{on } U_2\,,
  \end{cases}\\ \intertext{and}
  \abs{s'}^2_{\rho'} =
  \begin{cases}
    \log\, \rho'_1 & \text{on } U_1\,,\\
    \log \rho'_2 + \log\, \abs{s'_2}^2 & \text{on } U_2\,.
  \end{cases}\\
\end{gather*}
Furthermore, 
\begin{math}
  \log \, \abs{g'_{21}} = \log \, \abs{s'_2}
\end{math}
on $U_1\cap U_2$. 
Let us denote by $\norm{\cdot}=\abs{\cdot}_\rho= \abs{\cdot}_{\rho'}$
for simplicity. Using the relevant entries from~\eqref{eq:74} we have:
\begin{equation*}
  \begin{gathered}
    \theta^{0}_1 = c_1(\rho)\,\log \norm{s'}\,,\\
    \theta^{0}_2 = \delb\del \log \rho_2 \; \onehalf \log \rho'_2
    = dd^c \,\log\,\norm{s} \; \onehalf \log\rho'_2\,.    
  \end{gathered}
\end{equation*}
An elementary integration by parts leads to:
\begin{equation*}
  \begin{split}
    \int_X f_2 \theta^{0}_2
    &= \int_X f_2 \log\,\norm{s}\, c_1(\rho') \\
    &+ \int_X \log\,\norm{s} \, df_2 \wedge \onehalf d^c \log\rho'_2
    - \int_X \onehalf \log\rho'_2  \, df_2 \wedge d^c \log\,\norm{s}\,.
  \end{split}
\end{equation*}
On the other hand, we have
\begin{equation*}
  \int_X df_2 \wedge \theta^{1}_{12} =
  \int_X \log\,\norm{s}\; df_2\wedge d^c \log\,\abs{s'_2}
  -\int_X \log\,\abs{s'_2} \; df_2\wedge d^c \log\,\norm{s}\,,
\end{equation*}
and putting all terms together we obtain
\begin{equation}
  \label{eq:79}
  \begin{split}
    \int_X f_1\theta^{0}_1 +f_2\theta^{1}_2
    + df_2 \wedge \theta^{1}_{21}
    &=\int_X f_1\,c_1(\rho)\,\log \norm{s'}
    +\int_X f_2 \log\,\norm{s}\, c_1(\rho') \\
    &+ \int_X \log\,\norm{s}\; df_2\wedge d^c \log\,\norm{s'}
  +\int_X \log\,\norm{s'} \; df_1\wedge d^c \log\,\norm{s}
  \end{split}
\end{equation}
which, if expressed in terms of the squares of the norms, is (up
to a factor) the logarithm of $\norm{\dual{s}{s'}}$, as it is
found in \cite[formula 6.5.1]{MR89b:32038}. This version is due
to O. Gabber. Via the Poincar\'e-Lelong lemma (see,
e.g.\ \cite{gh:alg_geom})
\begin{equation*}
  \delb\del \log\,\norm{s} = c_1(\rho) + \tate\,[D]\,,
\end{equation*}
where $[D]$ is the delta-current supported at the divisor of $s$,
and similarly for $s'$, formula~\eqref{eq:79} can be recast into:
\begin{equation}
  \label{eq:80}
  \begin{split}
    \int_X f_1\theta^{0}_1 +f_2\theta^{1}_2
    + df_2 \wedge \theta^{1}_{21} &=
    \int_X dd^c \log\,\norm{s}\; \log\,\norm{s'}
    -\tate\,\log\,\norm{s} [D']
    -\tate\,\log\,\norm{s'} [D] \\
    &= -\tate\; \log\,\norm{\dual{s}{s'}}
  \end{split}
\end{equation}
which is what we wanted to show.
\end{proof}
This allows us to recast the Liouville functional for conformal metrics
on $X$ in the following form.
\begin{corollary}
  \label{corollary:2}
  The exponential of the Liouville functional defines an hermitian
  metric on the complex line $\dual{T_X}{T_X}$, namely for a conformal
  metric $\rho\in \cm{X}$ we have:
  \begin{equation*}
    \exp S[\log\rho] = \norm{\dual{T_X}{T_X}}\, \exp
    \tfrac{1}{2\pi} A_X(\rho)\,,
  \end{equation*}
  where $A_X(\rho)$ is the area of $X$ with respect to $\rho$.
\end{corollary}
\begin{remark}
  The above corollary justifies the choice made in
  Definition~\ref{def:2} for the various factors $\tate$.
\end{remark}
A similar result has been obtained in ref.\ \cite{zograf1990} by
considering the Liouville action functional defined on the Schottky
space, and in fact the results in loc.\ cit. are formulated in terms of
a Schottky \emph{family.}

Indeed, the statement in Corollary~\ref{corollary:2} can be immediately
reformulated for a family $\pi \colon X\rightarrow S$ with base
parameter space $S$ by considering the relative holomorphic tangent line
bundle $T_{X/S}$ with an hermitian fiber metric $\rho$. (Thus $\rho_s\in
\cm{X_s}$ for every fiber $X_s\,, s\in S$.) Notice that the fiber metric
$\rho$ needs not be critical (i.e.\ satisfying the fiberwise constant
negative curvature condition).

\section{Conclusions and outlook}
\label{sec:conclusions-outlook}

In this paper we have provided a geometric description of the
construction of the Liouville action functional for constant
negative curvature metrics on compact Riemann surfaces of genus
$g\geq 2$.

Our approach was to construct a Hermitian holomorphic Deligne
cohomology class as the cup square of the metrized holomorphic
tangent bundle $T_X$, and then show at the level of cocycles that
(modulo an area term) this construction agrees with those
in~\cite{math.CV/0204318} and earlier works.

Furthermore, our construction leads to the identification of the
class corresponding to the Liouville action with the determinant
of cohomology construction. Hence it could serve as an
alternative construction of the latter in terms of different
choices of cocycles.

One of the most important properties of the Liouville functional
from the works~\cite{zogtak1987-1,zogtak1987-2,zograf1990} is the
link with the Weil-Petersson geometry of the Teichm\"uller space.
From this point of view, an analysis of the behavior of the
Liouville action for families of Riemann surfaces is crucial. In
particular, a delicate analytic computation of the variation of
the Action with respect to the moduli was carried out in
the afore mentioned works to establish the link with the
Weil-Petersson metric.

We have not pursued these matters in the present work limiting
ourselves to establish the existence of a class for a metrized
bundle on the base of a family $\pi\colon X\to S$. A more precise
analysis would require, not only a more explicit description of
the map $\int_\pi$ (also advocated in~\cite{brymcl:deg4_II}), but
a full computation of the Leray sequence associated to the family
$\pi$. This is instrumental in defining \emph{relative} Hermitian
holomorphic classes, and in analyzing their variation with
respect to base parameters.

We hope to pursue this and other problems related to the
extension of the present work to singular metrics in a different
publication.

\appendix

\section{Cones}
\label{sec:Cones}

In the main body of the paper we have used iterated cones to
define the hermitian holomorphic Deligne complexes. One technical
problem one has to face concerns the homotopy (graded)
commutativity of the modified \bei\ product defined in
eq.~\eqref{eq:5}. A problem arises because the factors in the
cones are cones themselves and therefore they have multiplication
structures which are graded commutative up to homotopy to begin
with. We want to show that even in this situation the final
resulting product on cones is again homotopy graded commutative.
This ensures that on cohomology the product will be genuinely
graded commutative, so that in particular hermitian holomorphic
Deligne cohomology as defined in
section~\ref{sec:herm-holom-deligne} has the correct product
structure.

\subsection{Cones and homotopies}
\label{sec:cones-homotopies}

We consider the following situation. For $i=1,2,3$ we have maps of
complexes: 
\begin{math}
  f_i \colon A^\bullet_i \rightarrow B^\bullet_i\,,
\end{math}
and for $i<j$ maps
\begin{math}
  a_{ji} \colon A^\bullet_{i} \rightarrow A^\bullet_j
\end{math}
and
\begin{math}
  b_{ji} \colon B^\bullet_{i} \rightarrow B^\bullet_j\,.
\end{math}
Also, let
\begin{math}
  C^\bullet(f_i) =\cone (f_i:A^\bullet_i\rightarrow B^\bullet_i) \,,
\end{math}
for $i=1,2,3$. 

First, consider the homotopy commutative diagram:
\begin{equation}
  \label{eq:81}\vcenter{%
  \xymatrix{%
    A^\bullet_j \ar[r]^{f_j} \ar[d]_{a_{ij}} & B^\bullet_j
    \ar[d]^{b_{ij}} \ar@{=>}[dl]^{s_{ij}} \\
    A^\bullet_i \ar[r]_{f_i} & B^\bullet_i}}
\end{equation}
where
\begin{math}
  s_{ij} \colon A^\bullet_j \rightarrow B^{\bullet -1}_i
\end{math}
is the homotopy map of complexes:
\begin{equation*}
  f_ia_{ij} -b_{ij}f_j  = d\,s_{ij} + s_{ij}d\,.
\end{equation*}
An immediate verification yields
\begin{lemma}
  \label{lemma:cone-map}
  The diagram~\eqref{eq:81} can be extended to
  \begin{equation*}
    \xymatrix{%
      A^\bullet_j \ar[r]^{f_j} \ar[d]_{a_{ij}} &
      B^\bullet_j \ar[r] \ar[d]^{b_{ij}} \ar@{=>}[dl]^{s_{ij}} &
      C^\bullet (f_j) \ar[d]^{c_{ij}} \ar[r]^{[1]} &
      A^\bullet_j \ar[d]^{a_{ij}}\\
      A^\bullet_i \ar[r]_{f_i} & B^\bullet_i \ar[r] &
      C^\bullet (f_i) \ar[r]_{[1]} & A^\bullet_i
    }
  \end{equation*}
  where the map $c_{ij}$ is given by
  \begin{math}
    \bigl(
    \begin{smallmatrix}
      a_{ij} & 0 \\ -s_{ij} & b_{ij}
    \end{smallmatrix}
    \bigr)
  \end{math}
  and the squares containing the cones are in fact strictly commutative.
\end{lemma}
\begin{remark}
  This lemma is nothing other than the statement that any
  homotopy commutative diagram of the form~\eqref{eq:81} in the
  category of complexes in an abelian category can be extended to
  a (homotopy) commutative diagram of distinguished triangles,
  that is, one of the axioms defining a triangulated category,
  see, e.g. \cite{weibel_hom_alg}.
\end{remark}
For $k<j<i$ consider the homotopy commutative triangle
\begin{equation*}
  \xymatrix{%
    & A^\bullet_j  \ar[dr]^{a_{ij}} \ar@{=>}[d]^(.6){\alpha_{ijk}}& \\
    A^\bullet_k \ar[0,2]_{a_{ik}} \ar[ur]^{a_{jk}} & &
    A^\bullet_i}
\end{equation*}
where
\begin{math}
  a_{ik} -a_{ij}a_{jk} = d\alpha_{ijk} + \alpha_{ijk} d\,,
\end{math}
and similarly for the complexes $B^\bullet_i$, with a
corresponding homotopy $\beta_{ijk}$. Thus
\begin{math}
  \alpha_{ijk} \colon
  A^\bullet_k \rightarrow A^{\bullet -1}_i
\end{math}
and
\begin{math}
  \beta_{ijk} \colon
  B^\bullet_k \rightarrow B^{\bullet -1}_i\,.
\end{math}
Now consider the diagram:
\begin{equation}
  \label{eq:82}\vcenter{%
  \xymatrix@!{%
    A^\bullet_k \ar[0,2]^{a_{jk}} \ar[dr]_(.7){a_{ik}} \ar[2,0]_{f_k} & &
    A^\bullet_j \ar[dl]^{a_{ij}} \ar[2,0]^{f_j}
    \ar@{=>}[];[dl]+/ul 5ex/\\
    & A^\bullet_i \ar[2,0]^(.4){f_i} & \\
    B^\bullet_k \ar'[r]^(.7){b_{jk}}[rr] \ar[dr]_{b_{jk}}
    \ar@{=>}[ur] \ar@2{-->}[-1,0];[-2,1] & &
    B^\bullet_j \ar[dl]^{b_{ij}} \ar@{=>}[ul]
    \ar@{=>}[];[dl]+/ul 5ex/ \\
    & B^\bullet_i &
  }}
\end{equation}
The faces in~\eqref{eq:82} are homotopy commutative, however we assume
that composing the \emph{faces} is strictly commutative, namely the two
possible homotopies
\begin{math}
  b_{ij}\,b_{jk}\,f_k \Longrightarrow f_i\,a_{ik}
\end{math}
must be equal. Concretely, this corresponds to the relation
\begin{equation}
  \label{eq:83}
  s_{ik} + \beta_{ijk}\,f_{k} = f_i\,\alpha_{ijk}
  +s_{ij}\,a_{jk} +b_{ij}\,s_{jk}\,.
\end{equation}
We have:
\begin{lemma}
  \label{lemma:cone-homotopy}
  The map 
  \begin{equation*}
    \bigl(
    \begin{smallmatrix}
      -\alpha_{ijk} & 0 \\ 0 & \beta_{ijk}
    \end{smallmatrix}
    \bigr) \colon C^\bullet (f_k) \longrightarrow
    C^{\bullet -1}(f_i)
  \end{equation*}
  realizes the homotopy
  \begin{equation*}
    \xymatrix{%
      & C^\bullet (f_j) \ar[dr]^{c_{ij}} \ar@{=>}[d]& \\
      C^\bullet (f_k) \ar[ur]^{c_{jk}} \ar[rr]_{c_{ik}} & &
      C^\bullet (f_i)}
  \end{equation*}
\end{lemma}
\begin{proof}
  It is an elementary calculation based on writing 
  \begin{math}
    c_{ik} - c_{ik}\,c_{jk}
  \end{math}
  explicitly via the matrix representation given in
  Lemma~\ref{lemma:cone-map} and using eq.~\eqref{eq:83}.
\end{proof}

\subsection{Applications}
\label{sec:applications}

Consider the same setup as in section~\ref{sec:cup-products-cones}, with
the same complexes $X^\bullet_i$, etc., and diagrams:
\begin{equation*}
  \mathcal{D}_i \eqdef
  X^\bullet_i \overset{f_i}{\longrightarrow} Z^\bullet_i
  \overset{g_i}{\longleftarrow} Y^\bullet_i
\end{equation*}
from which we construct the cones
\begin{equation*}
  C(\mathcal{D}_i) =
  \cone (X^\bullet_i\oplus Y^\bullet_i
  \xrightarrow{f_i-g_i} Z^\bullet_i)[-1]\,,
  \quad i=1,2,3\,.
\end{equation*}
Moreover, following ref.\ \cite{bei:hodge_coho}, define
$\mathcal{D}_i\otimes \mathcal{D}_j$ by taking the tensor product
component-wise. Thus
\begin{equation*}
  \mathcal{D}_1\otimes \mathcal{D}_2 = 
  X^\bullet_1 \otimes X^\bullet_2
  \xrightarrow{f_1\otimes f_2}
  Z^\bullet_1\otimes Z^\bullet_2
  \xleftarrow{g_1\otimes g_2}
  Y^\bullet_1\otimes Y^\bullet_2\,.
\end{equation*}
Assuming as in section~\ref{sec:cup-products-cones} that the product
maps are compatible with the $f_i$, etc., the diagram
\begin{equation*}
  \mathcal{D}_1\otimes \mathcal{D}_2 \rightarrow \mathcal{D}_3
\end{equation*}
is of the same type as~\eqref{eq:81}, and therefore
lemma~\ref{lemma:cone-map} implies lemma~\ref{lemma:mod-cup-prod}.

Now, let the multiplication maps
\begin{math}
  X^\bullet_1 \otimes X^\bullet_2 \rightarrow X^\bullet_3
\end{math}
be graded commutative up to homotopy and similarly for the
$Y^\bullet_i$ and the $Z^\bullet_i$. We are interested in the
commutativity properties of multiplication map given by the \bei\
product~\eqref{eq:5}. 
\begin{proposition}
  \label{prop:4}
  The multiplication map
  \begin{math}
    \cup_\alpha\colon C(\mathcal{D}_1)\otimes C(\mathcal{D}_2)
    \longrightarrow C(\mathcal{D}_3)
  \end{math}
  given by~\eqref{eq:5} is homotopy graded commutative.
\end{proposition}
\begin{proof}
  The permutation operation on tensor products induces the diagram
  \begin{equation*}
    \xymatrix{%
      {\mathcal{D}_1}\otimes {\mathcal{D}_2}  \ar[dr] \ar[rr] & &
      {\mathcal{D}_2}\otimes {\mathcal{D}_1} \ar[dl] \\
      & {\mathcal{D}}_3 & }
  \end{equation*}
  which is of type~\eqref{eq:82} and we can apply
  lemma~\ref{lemma:cone-homotopy}.
\end{proof}
It follows from the proposition that the cohomology inherits a well
defined graded commutative product. This in particular applies to the
definition of hermitian holomorphic Deligne cohomology that uses the
cone~\eqref{eq:22}. Therefore we conclude that the cup
product~\eqref{eq:27} is graded commutative, as wanted.

\section{Conventions on Kleinian groups and fractional linear
  transformations} 
\label{sec:kleinian-groups}

As a reference the reader can consult, among many others, the
book \cite{MR99g:30055}. Let $\Gamma$ be a finitely generated
purely loxodromic non-elementary Kleinian group of the second
kind, so there is a nonempty discontinuity region
$U_\Gamma\subset \pione$. The limit set is $L_\Gamma = \pione
\setminus U_\Gamma$. According to Ahlfors' finiteness theorem the
quotient $X=U_\Gamma/\Gamma$ is a finite union of analytically
finite Riemann surfaces. Thus:
\begin{equation*}
  U_\Gamma/\Gamma = U_1/\Gamma_1 \sqcup \dots \sqcup U_n/\Gamma_n\,,
\end{equation*}
where $U_1,\dots,U_n$ are the inequivalent components of
$U_\Gamma$ and $\Gamma_1, \dots, \Gamma_n$ their stabilizers. By
way of example, a Schottky group has just one component, whereas
a Fuchsian or quasi-Fuchsian group will have exactly two
components.

We consider the map $U_\Gamma \to X$ as an \'etale cover, and
\cech\ cohomology with respect to it translates into group
cohomology for the group $\Gamma$, where the coefficient modules
are sections over $U_\Gamma$ of the relevant sheaves. The group
action is by pull-back.  According to the conventions of ref.\ 
\cite{aldtak2000} we will write the coboundary operation as:
\begin{equation}
  \label{eq:58}
    (\deltacheck c)_{\gamma_1,\dots,\gamma_n} =
    c_{\gamma_2,\dots,\gamma_n} 
    +\sum_{i=1}^{n-1}(-1)^i
    c_{\gamma_1,\dots,\gamma_i\gamma_{i+1},\dots,\gamma_n}
    +(-1)^n (c_{\gamma_1,\dots,\gamma_{n-1}})\cdot {\gamma_n}\,,
\end{equation}
where $c$ is an $n$-cochain with values in some \emph{right}
$\Gamma$-module $A$. The expression~\eqref{eq:58} is the \cech\ 
coboundary applied to the nerve of the cover $U_\Gamma \to X$.

We assume that $\Gamma$ is normalized, namely the point $\infty$
belongs to the limit set $L_\Gamma$. If $\gamma\in \Gamma$
corresponds to the fractional linear transformation:
\begin{equation*}
  \pione \ni z \longmapsto \gamma (z) = \frac{az +b}{cz +d}\,,
\end{equation*}
we have
\begin{equation}
  \label{eq:59}
  \gamma' (z) = \frac{\det \gamma}{c^2(z -z_\gamma)^2}\,,
\end{equation}
where $z_\gamma = -\frac{d}{c} \equiv \gamma^{-1}(\infty)$. Set
\begin{equation*}
  c(\gamma) \equiv c_\gamma \eqdef \frac{\det \gamma}{c^2}\,.
\end{equation*}
The following properties are easily verified. If $z_0$ and $z_\infty$
are the attracting and repelling fixed points, respectively, then
\begin{equation*}
  c_\gamma =
  \frac{(z_0 -z_\infty)^2\lambda_\gamma}{(1-\lambda_\gamma)^2} \,,
\end{equation*}
where $\lambda_\gamma$ is the dilating factor. For $\gamma_i$ and
$\gamma_j$ two elements of $\Gamma$, denote:
\begin{equation*}
  z_i = \gamma^{-1}(\infty)\,, \quad
  z_{ij} = (\gamma_i\gamma_j)^{-1}(\infty) =
  \gamma^{-1}_j (z_i)\,, \quad c_i = c_{\gamma_i}\,, \quad
  c_{ij} = c_{\gamma_i\gamma_j}\,.
\end{equation*}
We have the following relation:
\begin{equation}
  \label{eq:60}
  c_{\gamma_1\gamma_2}
  = \frac{c_{\gamma_1}}{c_{\gamma_2}}\, (z_{12} -z_2)^2\,.
\end{equation}
Finally, given four points $z_1,z_2, z_3, z_4\in \pione$, we define
their cross-ratio by:
\begin{equation}
  \label{eq:61}
  \CR{z_1}{z_2}{z_3}{z_4}
  = \frac{(z_1-z_2)(z_3-z_4)}{(z_1-z_4)(z_3-z_2)}\,.
\end{equation}

\bibliography{general}
\bibliographystyle{hamsplain}

\end{document}